\documentclass[ijoc,nonblindrev]{informs3}
\usepackage{mathrsfs}
\usepackage{latexsym}
\usepackage{amsmath}
\usepackage{theorem}
\usepackage{amssymb}
\usepackage{graphicx}
\usepackage{enumerate}
\usepackage{changebar}
\usepackage{url}
\usepackage{cases}
\usepackage{lscape}

\OneAndAHalfSpacedXII 

\usepackage{natbib}
 \bibpunct[, ]{(}{)}{,}{a}{}{,}%
 \def\bibfont{\small}%
\TheoremsNumberedThrough     

\EquationsNumberedThrough    

\newenvironment{preuve}{\noindent\textit{Proof. \ }}{\hfill$\Box$\medskip\par }
\newcounter{triang}
\newcommand{\toto}{ \refstepcounter{triang} (T\thetriang ) }
\newcounter{binary}

\newcommand{\ipopt}{\texttt{Ipopt}}
\newcommand{\csdp}{\texttt{csdp}}
\newcommand{\baron}{\texttt{Baron}}
\newcommand{\cplex}{\texttt{Cplex}}

\newcommand{\glomiqo}{\texttt{GloMIQO}}

\newcommand{\bb}{branch-and-bound}

\newcommand{\miqcr}{\texttt{MIQCR}}

\newcommand{\miqcrt}{\texttt{MIQCR-T}}

\newcommand{\setI}{\cal{I}}

\newcommand{\setSDP}{\cal{W}}
\newcommand{\setMC}{\cal{C}}
\newcommand{\setIMC}{\cal{U}}
\newcommand{\setGT}{\cal{G}}
\newcommand{\setT}{\cal{V}}
\newcommand{\setIndT}{T}

\newcommand{\setR}{\cal{R}}

\newcommand{\setC}{\cal{C}}

\newcommand{\OPT}{(OPT)}

\newcommand{\QP}{(P)}
\newcommand{\QPR}{(P_{S_0, \ldots, S_m})}
\newcommand{\QPT}{(P'_{S_0, \ldots, S_m})}
\newcommand{\QPTbar}{P'_{\bar{S}_0, \ldots , \bar{S}_m}}

\makeatletter
\let\markeq\ltx@label
\makeatother

\usepackage{url}

\begin{document}


\TITLE{Valid inequalities and global solution algorithm for Quadratically Constrained Quadratic Programs}

\RUNTITLE{Valid inequalities and global solution algorithm for QCQPs}
\ARTICLEAUTHORS{\AUTHOR{Am\'elie Lambert} \AFF{Cnam-CEDRIC, 292 Rue St Martin FR-75141 Paris Cedex 03, \EMAIL{amelie.lambert@cnam.fr}
}}

\RUNAUTHOR{Am\'elie Lambert}


\ABSTRACT{We consider the exact solution of problem $\QP$ that consists in minimizing a quadratic function subject to quadratic constraints. Starting from the classical convex relaxation that uses  the McCormick's envelopes, we introduce 12 inequalities that are derived from the ranges of the variables of $\QP$.  We prove that these general Triangle inequalities cut feasible solutions of the McCormick's envelopes. We then show how we can compute a convex relaxation $(P^*)$ which optimal value equals to the "Shor's plus RLT plus Triangle" semi-definite relaxation of $\QP$ that includes the new inequalities. We also propose a heuristic for solving this huge semi-definite program that serves as a separation algorithm. We then solve $\QP$ to global optimality with a \bb~based on  $(P^*)$. Moreover, as the new inequalities involved the lower and upper bounds on the original variables of $\QP$, their use in a \bb~framework accelerates the whole process. We show on the \texttt{unitbox} instances that our method outperforms the compared solvers.
}


\maketitle

\section{Introduction}

\par  In this manuscript, we aim at the exact solution of Mixed-Integer Quadratically Constrained Programs. These are the class of optimization problems where the objective function to minimize and the constraints are all quadratic. Such a problem can be formulated as follows:
\begin{numcases}{\QP} 
\min  f_0(x) \equiv \langle Q_0, xx^T \rangle   + c_0^Tx\nonumber\\
 \mbox{s.t.} \,\,\,\, f_r(x)  \equiv \langle Q_r, xx^T \rangle   + c_r^Tx\leq b_r & $ r \in \setR$ \label{cont}     \\
\qquad \ell_i \leq x_i \leq u_i & $ i \in \setI$ \label{bound}   \\
\qquad  x_i \in \mathbb{R} & $ i \in \setI$ \label{real}
\end{numcases}

\noindent  with $\langle A , B \rangle =\displaystyle{\sum_{i=1}^n}\displaystyle{\sum_{j=1}^n} a_{ij}b_{ij}$, and where $\setI=\{1,\ldots,n\}$, $\setR=\{1,\ldots,m\}$,  $\forall \,\, r \in \{0\} \cup \setR$, $(Q_r,c_r) \in \mathcal{S}_n \times \mathbb{R}^n$, with $\mathcal{S}_n$ the set of $n \times n$ symmetric matrices, $b \in \mathbb{R}^m$, and $u \in \mathbb{R}^n$.   Without loss of generality we suppose that $\ell \in \mathbb{R}_+^n$ and that  the feasible domain of $\QP$ is non-empty. Problem  $\QP$ trivially contains the case where there are quadratic equalities, since an equality can be replaced by two inequalities. It also contains the case of linear constraints since a linear equality is a quadratic constraint with a zero quadratic part.

\par Problem $\QP$ is a fundamental problem in global optimization. It arises in many applications, including facility location, production planning, multiperiod tankage quality problems in refinery processes, circle packing problems, euclidean distance geometry, or triangulation problems, see for instance~\cite{ATS99,AAT12,DorSah95,Pha82,LiuSah97,LocRab02,SutDai02,XieSah08}.

\par In the specific case where matrices $Q_r$  are positive semi-definite, $\QP$ is a convex problem that is polynomially solvable. In this case, efficient solvers are available. However, in general matrices $Q_r$ are indefinite and problem $\QP$ is $\mathcal{NP}$-hard~(\cite{GarJoh79}). In this case, the development of suitable relaxations is required for exact solution algorithms. Indeed, global optimization methods for solving $\QP$ are classically based on a \bb~framework in which a lower bound is computed by a certain relaxation scheme at each node of the search tree. The tightness of the relaxation and the efficiency in computing the lower bounds have a great impact on the behavior of such methods. Relaxation techniques for problems $\QP$ are mainly based on linearization, convex quadratic programming, or Semi-Definite Programming (SDP). Most of the proposed relaxations of the literature are either linear or quadratic and convex. To compute such a relaxation, the quadratic functions are reformulated as convex equivalent functions in an extended space of variables. More precisely, new variables $Y_{ij}$ are introduced for all $(i,j) \in \setI^2$, ($\setI^2$ is the cartesian product of a set $\setI$ by itself),  that are meant to satisfy the equalities $Y_{ij} = x_ix_j$.  The equivalent formulation in then solved by a \bb~algorithm based on a relaxation of the later non-convex equalities, for instance by linear constraints~(see for instance \cite{Cor76,AdaShe98,YajFuj98}). Software, implementing some of the methods described above, are available, see, for instance, \baron~(\cite{baron}), or \glomiqo~(\cite{MisFlo12,MisFlou2013,MSF15}). Using semi-definite relaxations within \bb~frameworks to solve $\QP$ was also widely studied~(\cite{Ans09,BurChe12,BurVan08,BurVan09,VanNem05,VanNem05b}). A semi-definite relaxation of $\QP$ can be obtained by lifting $x$ to a symmetric matrix $X = xx^T$ where the later non-convex constraints are relaxed to $X - xx^T \succeq 0$ ($M \succeq 0$ means that $M$ is positive semidefinite). This standard semi-definite relaxation is often referred to as "Shor's" relaxation of $\QP$. In~\cite{Ans09},  the "Shor's plus RLT" relaxation was introduced, where the McCormick's envelopes where added to the latter relaxation.  Method \miqcr~(\texttt{Mixed Integer Quadratic Convex Reformulation})~(\cite{EllLam19}) also handles $\QP$. In this approach, a tight quadratic convex relaxation to $\QP$ is calculated using the "Shor's plus RLT" relaxation. The original problem is then solved by a \bb~based on the obtained relaxation. 

\par In all these relaxations, the main common feature is that equalities  $Y_{ij} = x_ix_j$ are relaxed, and then forced by a \bb~algorithm in order to come back to the original problem.  The contributions of this paper are in order. In Section~\ref{subsec:validineq}, we start by the design of new families of valid inequalities for $\QP$. As the McCormick's envelopes, they are derived from the ranges $[\ell_i,u_i]$ of each variable $x_i$. After a complete description of these inequalities, we prove which of them cut feasible points of the McCormick's envelopes. We then link them to the literature, showing that in the specific case were $x \in [0,1]^n$, they amount to the well known Triangle inequalities introduced in~\cite{Pad89}. We thus call them General Triangle inequalities. Then in Section~\ref{sec:bestrelax}, we use the General Triangles to build quadratic convex relaxations of $\QP$ sharper than the ones used in method \miqcr. We also prove that we can compute a quadratic convex relaxation of $\QP$ that has the same optimal value as the "Shor's plus RLT plus Triangle" relaxation. Moreover, as the general upper and lower bounds $\ell$ and $u$ are involved within the new inequalities, the relaxation is again tighten in the course of the \bb~accelating the satisfaction of equalities $Y_{ij} = x_ix_j$. Since, there exists a huge number of new inequalities, we propose in Section~\ref{sec:bundle} to separate them within a  bundle algorithm that heuristically solves the  "Shor's plus RLT plus Triangle" relaxation. Finally in Section~\ref{sec:exp}, we evaluate the new method, called \miqcrt, on the instances of the literature and compare it to several solvers. Section~\ref{sec:conc} draws a conclusion.

\section{The General Triangle inequalities}
\label{subsec:validineq}

We start by building a convex relaxation of $\QP$ in an extended space of variables. As classically done, we introduce $\frac{n(n+1)}{2}$ new variables $Y_{ij}$ for all $(i,j) \in \setI^2$ that represent product $x_ix_j$. Then, we define sets $\setIMC=\{ (i,j) \in \setI^2 \, : \, i\leq j \}$ and $\setIndT = \{(i,j,t): (i,j) \in \setIMC, t=1, \ldots, 4 \}$, and we recall McCormick's envelopes $\setC = \Big \{(x,Y) \, : \, h^t_{ij}(x,Y) \leq 0\quad  \forall (i,j,t) \in \setIndT  \Big \} $ with:
\begin{numcases}{h^t_{ij}(x,Y) }
Y_{ij}-  u_jx_i - \ell_ix_j + u_j\ell_i  & $t = 1$ \label{mc1}\\
Y_{ij}-  u_ix_j - \ell_jx_i + u_i\ell_j  & $t = 2$\label{mc2}\\
-Y_{ij} + u_j x_i + u_i x_j-u_iu_j  & $t = 3$ \label{mc3}\\
- Y_{ij}  + \ell_jx_i + \ell_ix_j - \ell_i\ell_j  & $t=4$  \label{mc4}
 \end{numcases}
 We also consider any set of positive semi-definite matrices $S_0, \ldots ,S_m$. Then, $\forall \,\, r \in\{0\} \cup  \setR $, we formulate $f_r(x)$ as a sum of a quadratic function of the $x$ variables and a linear function of the $Y$ variables:
 $$f_{r,S_r}(x,Y) = \langle S_r, xx^T \rangle +c_r^Tx + \langle Q_r - S_r , Y \rangle .$$
It holds that $f_{r,S_r}(x,Y)$ is equal to $f_r(x)$ if $Y_{ij}=x_ix_j$. By replacing the initial functions $f_r(x)$ by the convex functions $f_{r,S_r}(x,Y)$, and by relaxing the non-convex equalities $Y_{ij} = x_ix_j$ with the inequalities of set $\setMC$, we obtain a family of quadratic convex relaxation of $\QP$:
\begin{numcases}{\QPR}
 \min \;\;\;  \langle S_0, xx^T \rangle +c_0^Tx + \langle Q_0 - S_0 , Y \rangle \nonumber \\
 \mbox{s.t.} \,\,\,\,   (\ref{bound})(\ref{real})  \nonumber \\
 \qquad  \langle S_r, xx^T \rangle +c_r^Tx + \langle Q_r - S_r , Y \rangle  \leq b_r  & $r \in \setR$ \nonumber \\
 \qquad h^t_{ij}(x,Y) \leq 0 & $(i,j,t) \in \setIndT $\nonumber \\
 \qquad Y_{jj} = Y_{ij} & $ (i,j) \in \overline{\setIMC}$ \nonumber
\end{numcases}
where $\overline{\setIMC}=\{ (i,j) \in \setI^2 \, : \, i < j \}$. Problem $\QPR$ is parameterized by the set of matrices $S_0, \ldots , S_m$, we observe that taking $\forall \,\, r \in\{0\} \cup  \setR $, $S_r=\mathbf{0}_n$  the zero $n \times n$ matrices amounts to the standard linearization of  $\QP$.
\bigskip

\par We now present new families of valid inequalities for $\QPR$ that strenghen this relaxation. As the McCormick's envelopes, they are derived from the ranges $[\ell_i,u_i]$ of each variable $x_i$.  The idea is to consider $\forall (i,j,k) \in \setT=\{ (i,j,k) \in I^3 \, : \, i<j<k \}$, three variables $x_i$, $x_j$ and $x_k$. Since these variables satisfy Constraints~(\ref{bound}), we have  $(u_i -x_i)(u_j - x_j)(u_k - x_k) \geq 0$, or equivalently:
$$ u_kx_ix_j + u_jx_ix_k + u_ix_jx_k -u_iu_kx_j - u_ju_kx_i - u_iu_jx_k + u_iu_ju_k  \geq x_ix_jx_k $$
using the McCormick inequality $x_jx_k \geq \ell_jx_k + \ell_kx_j -\ell_j\ell_k$, we get:
$$u_kx_ix_j + u_jx_ix_k + u_ix_jx_k -u_iu_kx_j - u_ju_kx_i - u_iu_jx_k + u_iu_ju_k \geq x_i(\ell_jx_k + \ell_kx_j -\ell_j\ell_k) $$
or equivalently the new quadratic inequality:
$$\boxed{ (\ell_k -u_k)x_ix_j + (\ell_j - u_j)x_ix_k - u_ix_jx_k + u_iu_kx_j + (u_ju_k - \ell_j\ell_k)x_i + u_iu_jx_k - u_iu_ju_k \leq 0}$$
that can be linearized by use of the variables $Y$. These inequalities are obviously valid by construction. In the example above, we also could have chosen to substitute the product $x_jx_k$ by its other McCormick envelope, i.e.  $x_jx_k \geq u_jx_k + u_kx_j -u_ju_k$,  or, to substitute either the product $x_ix_j$ or $x_ix_k$ by one of its two McCormick over estimators. Hence, considering all possible combinations for all $(i,j,k) \in \setT$ we obtain 8 families of 6 inequalities. The question is now to determine which of these 48 inequalities, when they are linearized, are non redundant in $\QPR$. We further prove in Propositions~\ref{prop1}--\ref{prop8} that 12 out of 48 of these inequalities cut feasible solutions of  $\QPR$.\\

\setcounter{equation}{0}
\noindent \textbf{ Family 1} We consider  $ (u_i -x_i)(u_j - x_j)(u_k - x_k) \geq 0$, and we get:\\
$\boxed{\Rightarrow   (\ell_k -u_k)Y_{ij} + (\ell_j - u_j)Y_{ik} - u_i Y_{jk}+ u_iu_kx_j + (u_ju_k - \ell_j\ell_k)x_i + u_iu_jx_k - u_iu_ju_k \leq 0 \quad \toto \label{c1a}}$\\
\noindent or symmetrically\\
$u_kx_ix_j + u_jx_ix_k + u_ix_jx_k -u_iu_kx_j - u_ju_kx_i - u_iu_jx_k + u_iu_ju_k  \geq x_ix_jx_k \geq x_j(\ell_ix_k + \ell_kx_i -\ell_i\ell_k)$\\
$\Rightarrow    (\ell_k -u_k)x_jx_i - u_jx_ix_k + ( \ell_i -u_i)x_jx_k + (u_iu_k - \ell_i\ell_k)x_j + u_ju_kx_i + u_iu_jx_k - u_iu_ju_k \leq 0 $\\
$\boxed{\Rightarrow    (\ell_k -u_k)Y_{ij} - u_jY_{ik}+ ( \ell_i -u_i)Y_{jk} + (u_iu_k - \ell_i\ell_k)x_j + u_ju_kx_i + u_iu_jx_k - u_iu_ju_k \leq 0 \quad \toto \label{c1b}}$\\
\noindent or  symmetrically\\
$u_kx_ix_j + u_jx_ix_k + u_ix_jx_k -u_iu_kx_j - u_ju_kx_i - u_iu_jx_k + u_iu_ju_k  \geq x_ix_jx_k \geq x_k(\ell_jx_i + \ell_ix_j -\ell_i\ell_j)$\\
$\Rightarrow   - u_kx_ix_j + (\ell_j - u_j)x_ix_k + ( \ell_i -u_i)x_jx_k + u_iu_kx_j + u_ju_kx_i + ( u_iu_j - \ell_i\ell_j)x_k - u_iu_ju_k \leq 0 $\\
$\boxed{\Rightarrow   - u_kY_{ij} + (\ell_j - u_j)Y_{ik}+ ( \ell_i -u_i)Y_{jk} + u_iu_kx_j + u_ju_kx_i + ( u_iu_j - \ell_i\ell_j)x_k - u_iu_ju_k \leq 0 \quad \toto \label{c1c}}$\\
\noindent or \\
$u_kx_ix_j + u_jx_ix_k + u_ix_jx_k -u_iu_kx_j - u_ju_kx_i - u_iu_jx_k + u_iu_ju_k  \geq x_ix_jx_k \geq x_i(u_jx_k + u_kx_j -u_ju_k)$\\
$\Rightarrow   u_ix_jx_k - u_iu_kx_j - u_iu_jx_k + u_iu_ju_k \geq 0 $\\
$\boxed{\Rightarrow   -u_iY_{jk} + u_iu_kx_j + u_iu_jx_k - u_iu_ju_k \leq 0 \quad \toto \label{c1d}}$\\
\noindent or symmetrically\\
$u_kx_ix_j + u_jx_ix_k + u_ix_jx_k -u_iu_kx_j - u_ju_kx_i - u_iu_jx_k + u_iu_ju_k  \geq x_ix_jx_k \geq x_j(u_ix_k + u_kx_i -u_iu_k)$\\
$\Rightarrow  u_jx_ix_k - u_ju_kx_i - u_iu_jx_k + u_iu_ju_k \geq 0 $\\
$\boxed{\Rightarrow   -u_jY_{ik} + u_ju_kx_i + u_iu_jx_k - u_iu_ju_k\leq 0 \quad \toto \label{c1e}}$\\
\noindent or symmetrically\\
$u_kx_ix_j + u_jx_ix_k + u_ix_jx_k -u_iu_kx_j - u_ju_kx_i - u_iu_jx_k + u_iu_ju_k  \geq x_ix_jx_k \geq x_k(u_ix_j + u_ix_j -u_iu_j)$\\
$\Rightarrow   u_kx_ix_j - u_ju_kx_i - u_iu_kx_j + u_iu_ju_k \geq 0 $\\
$\boxed{\Rightarrow   -u_kY_{ij} + u_ju_kx_i + u_iu_kx_j - u_iu_ju_k\leq 0 \quad \toto \label{c1f}}$\\

\begin{proposition}
  \label{prop1}
(T\ref{c1a})-(T\ref{c1c}) cut feasible solutions of $\QPR$, while~(T\ref{c1d})-(T\ref{c1f}) are redundant.
\end{proposition}
 
\begin{preuve}
\begin{enumerate}[(i)]
\item Consider the following solution  $(x,Y)$ satisfying $h^t_{ij}(x,Y) \leq 0$, for all $(i,j,t) \in T$:
\begin{itemize}
\item $x_i = \frac{u_i + \ell_i}{2}$, $x_j = \frac{u_j + \ell_j}{2}$, and $x_k = \frac{u_k + \ell_k}{2}$
\item $Y_{ij} = \frac{u_i\ell_j + \ell_iu_j}{2}$, $Y_{ik} = \frac{\ell_iu_k + u_i\ell_k}{2}$, and $Y_{jk} = \frac{u_j\ell_k + \ell_ju_k}{2}$
\end{itemize}
Using this solution in~(T\ref{c1a}), we get: $\frac{\ell_iu_j\ell_k}{2} + \frac{u_i\ell_j\ell_k}{2} - \frac{\ell_iu_ju_k}{2}  - \frac{u_i\ell_ju_k}{2} + \frac{\ell_i\ell_ju_k}{2} -  \frac{u_iu_j\ell_k}{2} +  \frac{u_iu_ju_k}{2}- \frac{\ell_i\ell_j\ell_k}{2}  = \frac{1}{2}(\ell_i - u_i)(u_j - \ell_j)(\ell_k - u_k) \geq 0$. This solution is thus cut off by the inequality~(T\ref{c1a}).  By symmetry, the proof is similar for~(T\ref{c1b}) and~(T\ref{c1c}).

\item (T\ref{c1d}) is equivalent to $ - u_i (Y_{jk} - u_kx_j - u_jx_k + u_ju_k) \leq 0$ that is redundant with~(\ref{mc3}). By symmetry, the proof is similar for~(T\ref{c1e}) and~(T\ref{c1f})

\end{enumerate}
\end{preuve}

\noindent \textbf{Family 2}  We consider $(u_i -x_i)(u_j - x_j)(x_k - \ell_k) \geq 0$, and we get:\\ 
\noindent  $u_jx_ix_k + \ell_kx_ix_j + u_ix_jx_k -u_j\ell_kx_i - u_i\ell_kx_j - u_iu_jx_k + u_iu_j\ell_k \leq x_ix_jx_k \leq x_i(u_kx_j + \ell_jx_k - \ell_ju_k)$ \\
$\boxed{\Rightarrow   (u_j - \ell_j) Y_{ik} + (\ell_k - u_k) Y_{ij} + u_iY_{jk}+ (\ell_ju_k -u_j\ell_k)x_i - u_i\ell_kx_j - u_iu_jx_k + u_iu_j\ell_k  \leq  0 \quad \toto \label{c2a}} $  
\noindent or  symmetrically\\
$u_jx_ix_k + \ell_kx_ix_j + u_ix_jx_k -u_j\ell_kx_i - u_i\ell_kx_j - u_iu_jx_k + u_iu_j\ell_k \leq x_ix_jx_k \leq x_j(u_kx_i + \ell_ix_k - u_k\ell_i)$ \\
$\boxed{\Rightarrow    u_j Y_{ik} + (\ell_k - u_k) Y_{ij} + (u_i-\ell_i)Y_{jk} -u_j\ell_kx_i + (\ell_iu_k - u_i\ell_k)x_j - u_iu_jx_k + u_iu_j\ell_k \leq  0\quad \toto \label{c2b}} $  

\noindent or\\
$u_jx_ix_k + \ell_kx_ix_j + u_ix_jx_k -u_j\ell_kx_i - u_i\ell_kx_j - u_iu_jx_k + u_iu_j\ell_k \leq x_ix_jx_k \leq x_k(u_ix_j + \ell_jx_i - u_i\ell_j)$ \\
$\Rightarrow   (u_j - \ell_j) x_ix_k + \ell_k x_ix_j -u_j\ell_kx_i - \ell_ku_ix_j + u_i(\ell_j- u_j)x_k + u_iu_j \ell_k  \leq  0 $ \\
$\boxed{\Rightarrow  (u_j - \ell_j) Y_{ik} + \ell_k Y_{ij} -u_j\ell_kx_i - u_i\ell_kx_j +u_i(\ell_j- u_j)x_k + u_iu_j \ell_k  \leq  0   \quad \toto \label{c2c}}$

\noindent or symmetrically \\
$u_jx_ix_k + \ell_kx_ix_j + u_ix_jx_k -u_j\ell_kx_i - u_i\ell_kx_j - u_iu_jx_k + u_iu_j\ell_k \leq x_ix_jx_k \leq x_k(u_jx_i + \ell_ix_j - \ell_iu_j)$ \\
$\Rightarrow    \ell_k x_ix_j + (u_i - \ell_i) x_jx_k  - u_j\ell_kx_i -u_i\ell_kx_j +(\ell_i- u_i)u_jx_k + u_iu_j \ell_k  \leq  0 $ \\
$\boxed{\Rightarrow    \ell_k Y_{ij} + (u_i - \ell_i) Y_{jk}  - u_j\ell_kx_i -u_i\ell_kx_j+(\ell_i- u_i)u_jx_k + u_iu_j \ell_k  \leq  0   \quad \toto \label{c2d}}$  

\noindent or \\
$ u_jx_ix_k + \ell_kx_ix_j + u_ix_jx_k -u_j\ell_kx_i - u_i\ell_kx_j - u_iu_jx_k + u_iu_j\ell_k \leq x_ix_jx_k \leq x_i(u_jx_k + \ell_kx_j - u_j\ell_k)$\\
$ \Rightarrow  u_ix_jx_k - u_i\ell_kx_j - u_iu_jx_k + u_iu_j\ell_k \leq 0$\\
$ \boxed{\Rightarrow  u_iY_{jk} - u_i\ell_kx_j - u_iu_jx_k + u_iu_j\ell_k  \leq 0 \quad \toto \label{c2e}}$

\noindent or symmetrically\\
$ u_jx_ix_k + \ell_kx_ix_j + u_ix_jx_k -u_j\ell_kx_i - u_i\ell_kx_j - u_iu_jx_k + u_iu_j\ell_k \leq x_ix_jx_k \leq x_j(u_ix_k + \ell_kx_i - u_i\ell_k)$\\
$ \Rightarrow  u_jx_ix_k - u_j\ell_kx_i - u_iu_jx_k + u_iu_j\ell_k \leq 0$\\
$ \boxed{\Rightarrow  u_jY_{ik}  - u_j\ell_kx_i - u_iu_jx_k + u_iu_j\ell_k   \leq 0 \quad \toto \label{c2f}}$\\

\begin{proposition}
    \label{prop2}
(T\ref{c2a}) and~(T\ref{c2b}) cut feasible solutions of $\QPR$, while~(T\ref{c2c})-(T\ref{c2f}) are redundant.
\end{proposition}
 
\begin{preuve}

\begin{enumerate}[(i)]
\item Consider the following solution  $(x,Y)$ satisfaying $h^t_{ij}(x,Y) \leq 0$, for all $(i,j,t) \in T$:
\begin{itemize}
\item $x_i = \frac{u_i + \ell_i}{2}$, $x_j = \frac{u_j + \ell_j}{2}$, and $x_k = \frac{u_k + \ell_k}{2}$
\item $Y_{ij} = \frac{u_i\ell_j + \ell_iu_j}{2}$, $Y_{ik} = \frac{u_iu_k + \ell_i\ell_k}{2}$, and $Y_{jk} = \frac{u_ju_k + \ell_j\ell_k}{2}$
\end{itemize}
Using this solution in~(T\ref{c2a}), we get: $\frac{u_iu_ju_k}{2} + \frac{\ell_iu_j\ell_k}{2} - \frac{u_i\ell_ju_k}{2}  - \frac{\ell_i\ell_j\ell_k}{2} + \frac{u_i\ell_j\ell_k}{2} -  \frac{\ell_iu_ju_k}{2} +  \frac{\ell_i\ell_ju_k}{2}- \frac{u_iu_j\ell_k}{2} = \frac{1}{2}(\ell_i - u_i)(u_j - \ell_j)(\ell_k - u_k) \geq 0$. This solution is thus cut off by the inequality~(T\ref{c2a}). By symmetry, the proof is similar for~(T\ref{c2b}).w

\item (T\ref{c2c}) is equivalent to  $(u_j - \ell_j)(Y_{ik} - u_ix_k - \ell_kx_i + u_i\ell_k) +  \ell_k(Y_{ij}  - \ell_jx_i - u_ix_j+ u_i\ell_j) \leq 0$ that is redundant with~(\ref{mc2}).  By symmetry, the proof is similar for~(T\ref{c2d}).
\item  (T\ref{c2e}) is equivalent to $u_i (Y_{jk} - \ell_kx_j - u_jx_k + u_j\ell_k) \leq 0$ that is redundant with~(\ref{mc2}). By symmetry, the proof is similar for~(T\ref{c2f}).
\end{enumerate}
\end{preuve}

\noindent \textbf{Family 3} We consider $(u_i -x_i)(x_j - \ell_j)(u_k - x_k) \geq 0$ and we get:\\
\noindent  $u_kx_ix_j + u_ix_jx_k + \ell_jx_ix_k -u_iu_kx_j - \ell_ju_kx_i - u_i\ell_jx_k + u_i\ell_ju_k \leq x_ix_jx_k \leq  x_i(u_kx_j + \ell_jx_k - \ell_ju_k)$ \\
$ \Rightarrow u_ix_jx_k  -u_iu_kx_j  - u_i\ell_jx_k + u_i\ell_ju_k \leq 0$\\
$\boxed{\Rightarrow     u_iY_{jk}  -u_iu_kx_j  - u_i\ell_jx_k + u_i\ell_ju_k  \leq  0 \quad \toto \label{c3a}} $

\noindent or symmetrically\\
$u_kx_ix_j + u_ix_jx_k + \ell_jx_ix_k -u_iu_kx_j - \ell_ju_kx_i - u_i\ell_jx_k + u_i\ell_ju_k \leq x_ix_jx_k \leq  x_k(u_ix_j + \ell_jx_i - u_i\ell_j)$ \\
$\Rightarrow  u_kx_ix_j -u_iu_kx_j - \ell_ju_kx_i + u_i\ell_ju_k \leq  0$\\
$ \boxed{\Rightarrow    u_kY_{ij} - u_iu_kx_j - \ell_ju_kx_i + u_i\ell_ju_k  \leq  0 \quad \toto \label{c3b}}$

\noindent or  \\
$u_kx_ix_j + u_ix_jx_k + \ell_jx_ix_k -u_iu_kx_j - \ell_ju_kx_i - u_i\ell_jx_k + u_i\ell_ju_k \leq x_ix_jx_k\leq x_j(u_kx_i + \ell_ix_k - \ell_iu_k)$ \\
$\Rightarrow (u_i - \ell_i)x_jx_k + \ell_jx_ix_k + u_k(\ell_i -u_i)x_j - \ell_ju_kx_i - u_i\ell_jx_k + u_i\ell_ju_k \leq 0$\\
$\boxed{\Rightarrow   (u_i - \ell_i)Y_{jk} + \ell_jY_{ij} +u_k(\ell_i -u_i)x_j - \ell_ju_kx_i - u_i\ell_jx_k + u_i\ell_ju_k\leq 0 \quad \toto \label{c3c} }$

\noindent or  symmetrically\\
$u_kx_ix_j + u_ix_jx_k + \ell_jx_ix_k -u_iu_kx_j - \ell_ju_kx_i - u_i\ell_jx_k + u_i\ell_ju_k \leq x_ix_jx_k\leq x_j(u_ix_k + \ell_kx_i - u_i\ell_k)$ \\
$\Rightarrow  (u_k -\ell_k)x_ix_j + \ell_jx_ix_k +u_i(\ell_k -u_k)x_j - \ell_ju_kx_i - u_i\ell_jx_k + u_i\ell_ju_k \leq 0$\\
$\boxed{\Rightarrow (u_k -\ell_k)Y_{ij} + \ell_jY_{ik} +u_i(\ell_k -u_k)x_j - \ell_ju_kx_i - u_i\ell_jx_k + u_i\ell_ju_k \leq 0 \quad \toto \label{c3d} }$  

\noindent or  \\
$u_kx_ix_j + u_ix_jx_k + \ell_jx_ix_k -u_iu_kx_j - \ell_ju_kx_i - u_i\ell_jx_k + u_i\ell_ju_k \leq x_ix_jx_k \leq  x_i(u_jx_k + \ell_kx_j - u_j\ell_k)$ \\
$\boxed{(u_k - \ell_k) Y_{ij} + u_iY_{jk} + (\ell_j - u_j) Y_{ik} -u_iu_kx_j +(u_j\ell_k - \ell_ju_k)x_i - u_i\ell_jx_k + u_i\ell_ju_k \leq 0 \quad \toto \label{c3e}}$

\noindent or  symmetrically\\
$u_kx_ix_j + u_ix_jx_k + \ell_jx_ix_k -u_iu_kx_j - \ell_ju_kx_i - u_i\ell_jx_k + u_i\ell_ju_k \leq x_ix_jx_k \leq  x_k(u_jx_i + \ell_ix_j - \ell_iu_j)$ \\
$\boxed{u_kY_{ij} + (u_i - \ell_i) Y_{jk} + (\ell_j - u_j) Y_{ik} -u_iu_kx_j - \ell_ju_kx_i +(\ell_iu_j - u_i\ell_j)x_k + u_i\ell_ju_k \leq 0 \quad \toto \label{c3f}}$\\

\begin{proposition}
    \label{prop3}
(T\ref{c3e}) and~(T\ref{c3f}) cut feasible solutions of $\QPR$, while~(T\ref{c3a})-(T\ref{c3d}) are redundant.
\end{proposition}
 
\begin{preuve}

\begin{enumerate}[(i)]
\item Consider the following solution  $(x,Y)$ satisfaying $h^t_{ij}(x,Y) \leq 0$, for all $(i,j,t) \in T$:
\begin{itemize}
\item $x_i = \frac{u_i + \ell_i}{2}$, $x_j = \frac{u_j + \ell_j}{2}$, and $x_k = \frac{u_k + \ell_k}{2}$
\item $Y_{ij} = \frac{\ell_i\ell_j + u_iu_j}{2}$, $Y_{ik} = \frac{u_i\ell_k + \ell_iu_k}{2}$, and $Y_{jk} = \frac{u_ju_k + \ell_j\ell_k}{2}$
\end{itemize}
Using this solution in Constraints~(T\ref{c3e}), we get: $\frac{u_iu_ju_k}{2} + \frac{\ell_iu_j\ell_k}{2} - \frac{u_i\ell_ju_k}{2}  - \frac{\ell_i\ell_j\ell_k}{2} + \frac{u_i\ell_j\ell_k}{2} -  \frac{\ell_iu_ju_k}{2} +  \frac{\ell_i\ell_ju_k}{2}- \frac{u_iu_j\ell_k}{2} = \frac{1}{2}(\ell_i - u_i)(u_j - \ell_j)(\ell_k - u_k) \geq 0$. This solution is thus cut off by the inequality~(T\ref{c3e}). By symmetry, the proof is similar for~(T\ref{c3f}).  
\item (T\ref{c3a}) is equivalent to $u_i (Y_{jk} - u_kx_j - \ell_jx_k + \ell_ju_k) \leq 0$ that is redundant with~(\ref{mc1}). By symmetry, the proof is similar for~(T\ref{c3b}). 
\item (T\ref{c3c}) is equivalent to $(u_i - \ell_i)(Y_{jk} - u_kx_j - \ell_jx_k + \ell_ju_k) +  \ell_j(Y_{ik}  - \ell_ix_k - u_kx_i+ \ell_iu_k) \leq 0$ that is redundant with~(\ref{mc1}). By symmetry, the proof is similar for~(T\ref{c3d}). 
\end{enumerate}
\end{preuve}

\noindent \textbf{ Family 4} We consider $(x_i -\ell_i)(u_j - x_j)(u_k - x_k) \geq 0$ and we get:\\
\noindent  $u_kx_ix_j + \ell_ix_jx_k + u_jx_ix_k - \ell_iu_kx_j - u_ju_kx_i - \ell_iu_jx_k + \ell_iu_ju_k \leq x_ix_jx_k \leq x_i(u_kx_j + \ell_jx_k - \ell_ju_k)$ \\
$\Rightarrow  \ell_ix_jx_k + (u_j - \ell_j)x_ix_k - \ell_iu_kx_j + u_k(\ell_j- u_j)x_i - \ell_iu_jx_k + \ell_iu_ju_k \leq 0$ \\
$\boxed{\Rightarrow  \ell_iY_{jk} + (u_j - \ell_j)Y_{ik} - \ell_iu_kx_j + u_k(\ell_j- u_j)x_i - \ell_iu_jx_k + \ell_iu_ju_k \leq 0  \quad \toto \label{c4a}}$  \\
\noindent or symmetrically\\
$u_kx_ix_j + \ell_ix_jx_k + u_jx_ix_k - \ell_iu_kx_j - u_ju_kx_i - \ell_iu_jx_k + \ell_iu_ju_k \leq x_ix_jx_k \leq x_i(u_jx_k + \ell_kx_j - u_j\ell_k)$ \\
$\Rightarrow (u_k - \ell_k)x_ix_j + \ell_ix_jx_k  - \ell_iu_kx_j +u_j(\ell_k - u_k)x_i - \ell_iu_jx_k + \ell_iu_ju_k \leq 0 $\\
$\boxed{\Rightarrow (u_k - \ell_k)Y_{ij} + \ell_iY_{jk}  - \ell_iu_kx_j +u_j(\ell_k - u_k)x_i - \ell_iu_jx_k + \ell_iu_ju_k \leq 0  \quad \toto \label{c4b}}$  \\
\noindent or \\
$u_kx_ix_j + \ell_ix_jx_k + u_jx_ix_k - \ell_iu_kx_j - u_ju_kx_i - \ell_iu_jx_k + \ell_iu_ju_k \leq x_ix_jx_k \leq x_j(u_kx_i + \ell_ix_k - \ell_iu_k)$ \\
$\Rightarrow u_jx_ix_k -  u_ju_kx_i - \ell_iu_jx_k + \ell_iu_ju_k \leq 0 $\\
$ \boxed{\Rightarrow  u_jY_{ik} -  u_ju_kx_i - \ell_iu_jx_k + \ell_iu_ju_k \leq 0\quad \toto \label{c4c} }$  \\
\noindent or symmetrically\\
$u_kx_ix_j + \ell_ix_jx_k + u_jx_ix_k - \ell_iu_kx_j - u_ju_kx_i - \ell_iu_jx_k + \ell_iu_ju_k \leq x_ix_jx_k \leq x_k(u_jx_i + \ell_ix_j - \ell_iu_j)$ \\
$\Rightarrow u_kx_ix_j - \ell_iu_kx_j - u_ju_kx_i  + \ell_iu_ju_k \leq 0$ \\
$ \boxed{\Rightarrow u_kY_{ij} - \ell_iu_kx_j - u_ju_kx_i  + \ell_iu_ju_k  \leq 0\quad \toto \label{c4d} }$  \\
\noindent or \\
$u_kx_ix_j + \ell_ix_jx_k + u_jx_ix_k - \ell_iu_kx_j - u_ju_kx_i - \ell_iu_jx_k + \ell_iu_ju_k \leq x_ix_jx_k \leq x_j(u_ix_k + \ell_kx_i - u_i\ell_k)$ \\
$ \boxed{\Rightarrow(u_k - \ell_k)Y_{ij} + (\ell_i - u_i) Y_{jk} + u_jY_{ik} +(u_i\ell_k - \ell_iu_k)x_j - u_ju_kx_i - \ell_iu_jx_k + \ell_iu_ju_k \leq 0 \quad \toto \label{c4e}}$ \\
\noindent or symmetrically\\
$u_kx_ix_j + \ell_ix_jx_k + u_jx_ix_k - \ell_iu_kx_j - u_ju_kx_i - \ell_iu_jx_k + \ell_iu_ju_k \leq x_ix_jx_k \leq x_k(u_ix_j + \ell_jx_i - u_i\ell_j)$ \\
$ \boxed{\Rightarrow u_kY_{ij} + (\ell_i - u_i) Y_{jk} + (u_j - \ell_j)Y_{ik} - \ell_iu_kx_j - u_ju_kx_i + (u_i\ell_j - \ell_iu_j) x_k + \ell_iu_ju_k  \leq 0 \quad \toto \label{c4f}}$

\begin{proposition}
    \label{prop4}
(T\ref{c4e}) and~(T\ref{c4f}) cut feasible solutions of $\QPR$, while~(T\ref{c4a})-(T\ref{c4d}) are redundant.
\end{proposition}
 
\begin{preuve}

\begin{enumerate}[(i)]
\item Consider the following solution  $(x,Y)$ satisfaying $h^t_{ij}(x,Y) \leq 0$, for all $(i,j,t) \in T$:
\begin{itemize}
\item $x_i = \frac{u_i + \ell_i}{2}$, $x_j = \frac{u_j + \ell_j}{2}$, and $x_k = \frac{u_k + \ell_k}{2}$
\item $Y_{ij} = \frac{\ell_i\ell_j + u_iu_j}{2}$, $Y_{ik} = \frac{\ell_i\ell_k + u_iu_k}{2}$, and $Y_{jk} = \frac{u_j\ell_k + \ell_ju_k}{2}$
\end{itemize}
Using this solution in Constraints~(T\ref{c4e}), we get: $\frac{u_iu_ju_k}{2} + \frac{\ell_iu_j\ell_k}{2} - \frac{u_i\ell_ju_k}{2}  - \frac{\ell_i\ell_j\ell_k}{2} + \frac{u_i\ell_j\ell_k}{2} -  \frac{\ell_iu_ju_k}{2} +  \frac{\ell_i\ell_ju_k}{2}- \frac{u_iu_j\ell_k}{2} = \frac{1}{2}(\ell_i - u_i) (\ell_j - u_j)(u_k - \ell_k)\geq 0$. This solution is thus cut off by the inequality~(T\ref{c4e}).  By symmetry, the proof is similar for~(T\ref{c4f}).
\item (T\ref{c4a}) is equivalent to $(u_j - \ell_j)(Y_{ik} - u_kx_i - \ell_ix_k + \ell_iu_k) +  \ell_i(Y_{jk}  - \ell_jx_k - u_kx_j+ \ell_ju_k) \leq 0$ that is redundant with~(\ref{mc1}). By symmetry, the proof is similar for~(T\ref{c4b}).
\item (T\ref{c4c}) is equivalent to $u_j (Y_{ik} - u_kx_i - \ell_ix_k + \ell_iu_k) \leq 0 $ that is redundant with~(\ref{mc1}). By symmetry, the proof is similar for~(T\ref{c4d}).
\end{enumerate}
\end{preuve}

\noindent \textbf{ Family 5} We consider $(u_i -x_i)(x_j - \ell_j)(x_k - \ell_k) \geq 0$ and we get:\\
\noindent  $u_ix_jx_k + \ell_jx_ix_k + \ell_kx_ix_j - u_i\ell_kx_j - \ell_j\ell_kx_i - u_i\ell_jx_k + u_i\ell_j\ell_k \geq x_ix_jx_k \geq x_i(\ell_jx_k + \ell_kx_j - \ell_j\ell_k)$ \\
$\Rightarrow - u_ix_jx_k + u_i\ell_kx_j + u_i\ell_jx_k - u_i\ell_j\ell_k \leq 0$ \\
$\boxed{\Rightarrow  - u_iY_{jk}+ u_i\ell_kx_j + u_i\ell_jx_k - u_i\ell_j\ell_k \leq 0 \quad  \toto \label{c5a}}$  \\
\noindent or \\
$u_ix_jx_k + \ell_jx_ix_k + \ell_kx_ix_j - u_i\ell_kx_j - \ell_j\ell_kx_i - u_i\ell_jx_k + u_i\ell_j\ell_k \geq x_ix_jx_k \geq x_j(u_ix_k + u_kx_i - u_iu_k)$ \\
$\Rightarrow - \ell_jx_ix_k + (u_k - \ell_k)x_ix_j + (\ell_k - u_k)u_ix_j + \ell_j\ell_kx_i + u_i\ell_jx_k - u_i\ell_j\ell_k \leq 0$ \\
$\boxed{\Rightarrow - \ell_jY_{ik} + (u_k - \ell_k)Y_{ij} + (\ell_k - u_k)u_ix_j + \ell_j\ell_kx_i + u_i\ell_jx_k - u_i\ell_j\ell_k \leq 0\quad  \toto \label{c5b} }$  \\
\noindent or symmetrically\\
$u_ix_jx_k + \ell_jx_ix_k + \ell_kx_ix_j - u_i\ell_kx_j - \ell_j\ell_kx_i - u_i\ell_jx_k + u_i\ell_j\ell_k \geq x_ix_jx_k \geq x_k(u_jx_i + u_ix_j - u_iu_j)$ \\
$\Rightarrow (u_j \ell_j)x_ix_k - \ell_kx_ix_j + u_i\ell_kx_j + \ell_j\ell_kx_i + u_i(\ell_j - u_j)x_k - u_i\ell_j\ell_k \leq 0 $\\
$\boxed{\Rightarrow (u_j \ell_j)Y_{ik} - \ell_kY_{ij} + u_i\ell_kx_j + \ell_j\ell_kx_i + u_i(\ell_j - u_j)x_k - u_i\ell_j\ell_k   \leq 0 \quad  \toto \label{c5c} }$  \\
\noindent or  symmetrically\\
$u_ix_jx_k + \ell_jx_ix_k + \ell_kx_ix_j - u_i\ell_kx_j - \ell_j\ell_kx_i - u_i\ell_jx_k + u_i\ell_j\ell_k \geq x_ix_jx_k \geq x_j(\ell_ix_k + \ell_kx_i - \ell_i\ell_k)$ \\
$\Rightarrow (\ell_i - u_i)x_jx_k - \ell_jx_ix_k + (u_i - \ell_i)\ell_kx_j + \ell_j\ell_kx_i + u_i\ell_jx_k - u_i\ell_j\ell_k \leq 0 $ \\
$\boxed{\Rightarrow  (\ell_i - u_i)Y_{jk} - \ell_jY_{ik} + (u_i - \ell_i)\ell_kx_j + \ell_j\ell_kx_i + u_i\ell_jx_k - u_i\ell_j\ell_k \leq 0 \quad  \toto \label{c5d} }$  \\
\noindent or  symmetrically\\
$u_ix_jx_k + \ell_jx_ix_k + \ell_kx_ix_j - u_i\ell_kx_j - \ell_j\ell_kx_i - u_i\ell_jx_k + u_i\ell_j\ell_k \geq x_ix_jx_k \geq x_k(\ell_ix_j + \ell_jx_i - \ell_i\ell_j)$ \\
$\Rightarrow (\ell_i - u_i)x_jx_k - \ell_kx_ix_j + u_i\ell_kx_j + \ell_j\ell_kx_i +(u_i -\ell_i)\ell_jx_k - u_i\ell_j\ell_k \leq 0$ \\
$ \boxed{\Rightarrow (\ell_i - u_i)Y_{jk} - \ell_kY_{ij} + u_i\ell_kx_j + \ell_j\ell_kx_i +(u_i -\ell_i)\ell_jx_k - u_i\ell_j\ell_k \leq 0 \quad  \toto \label{c5e} }$\\
\noindent or \\
$u_ix_jx_k + \ell_jx_ix_k + \ell_kx_ix_j - u_i\ell_kx_j - \ell_j\ell_kx_i - u_i\ell_jx_k + u_i\ell_j\ell_k   \geq x_ix_jx_k \geq x_i(u_jx_k + u_kx_j - u_ju_k)$ \\
$ \boxed{\Rightarrow  - u_i Y_{jk} + (u_j - \ell_j)Y_{ik} + (u_k - \ell_k) Y_{ij} + u_i\ell_kx_j + (\ell_j\ell_k -u_j u_k)x_i + u_i\ell_jx_k - u_i\ell_j\ell_k \leq 0\quad  \toto \label{c5f}}$ \\

\begin{proposition}
    \label{prop5}
(T\ref{c5f}) cuts feasible solutions of $\QPR$, while~(T\ref{c5a})-(T\ref{c5e}) are redundant.
\end{proposition}
 
\begin{preuve}

\begin{enumerate}[(i)]
\item Consider the following solution  $(x,Y)$ satisfaying $h^t_{ij}(x,Y) \leq 0$, for all $(i,j,t) \in T$:
\begin{itemize}
\item $x_i = \frac{u_i + \ell_i}{2}$, $x_j = \frac{u_j + \ell_j}{2}$, and $x_k = \frac{u_k + \ell_k}{2}$
\item $Y_{ij} = \frac{u_iu_j + \ell_i\ell_j}{2}$, $Y_{ik} = \frac{\ell_i\ell_k + u_iu_k}{2}$, and $Y_{jk} = \frac{u_j\ell_k + \ell_ju_k}{2}$
\end{itemize}
Using this solution in Constraints~(T\ref{c5f}), we get: $\frac{\ell_iu_j\ell_k}{2} - \frac{u_i\ell_ju_k}{2} + \frac{u_iu_ju_k}{2}  + \frac{\ell_i\ell_ju_k}{2} - \frac{u_iu_j\ell_k}{2} -  \frac{\ell_i\ell_j\ell_k}{2} +  \frac{\ell_iu_ju_k}{2} + \frac{u_i\ell_j\ell_k}{2} = \frac{1}{2}(\ell_i - u_i)(u_j - \ell_j)(\ell_k - u_k) \geq 0$. This solution is thus cut off by the inequality~(T\ref{c5f}). 
\item  (T\ref{c5a}) is  equivalent to $u_i (- Y_{jk} + \ell_kx_j + \ell_jx_k - \ell_j\ell_k) \leq 0 $ that is redundant with~(\ref{mc4}).
\item  (T\ref{c5b}) is equivalent to $(u_k - \ell_k)(Y_{ij} - u_ix_j - \ell_jx_i + u_i\ell_j) +  \ell_j(- Y_{ik}  + u_ix_k + u_kx_i+ u_iu_k) \leq 0$ that is redundant with~(\ref{mc2}) and~(\ref{mc3}). By symmetry, the proof is similar for~(T\ref{c5c})-(T\ref{c5e}).
\end{enumerate}
\end{preuve}

\noindent \textbf{Family 6}  We consider  $(x_i -\ell_i)(u_j - x_j)(x_k - \ell_k) \geq 0$ and we get:\\
\noindent  $\ell_ix_jx_k + u_jx_ix_k + \ell_kx_ix_j - \ell_i\ell_kx_j - u_j\ell_kx_i - \ell_iu_jx_k + \ell_iu_j\ell_k \geq x_ix_jx_k \geq x_i(\ell_jx_k + \ell_kx_j - \ell_j\ell_k)$ \\
$\Rightarrow  -\ell_ix_jx_k + (\ell_j - u_j)x_ix_k + \ell_i\ell_kx_j + (u_j - \ell_j)\ell_kx_i + \ell_iu_jx_k - \ell_iu_j\ell_k \leq 0 $ \\
$\boxed{\Rightarrow  -\ell_iY_{jk} + (\ell_j - u_j)Y_{ik} + \ell_i\ell_kx_j + (u_j - \ell_j)\ell_kx_i + \ell_iu_jx_k - \ell_iu_j\ell_k \leq 0 \toto \label{c6a}}$  \\
\noindent or symmetrically\\
$\ell_ix_jx_k + u_jx_ix_k + \ell_kx_ix_j - \ell_i\ell_kx_j - u_j\ell_kx_i - \ell_iu_jx_k + \ell_iu_j\ell_k \geq x_ix_jx_k \geq x_k(\ell_jx_i + \ell_ix_j - \ell_i\ell_j)$ \\
$\Rightarrow  (\ell_j - u_j)x_ix_k - \ell_kx_ix_j + \ell_i\ell_kx_j + u_j\ell_kx_i +\ell_i(u_j - \ell_j)x_k - \ell_iu_j\ell_k \leq 0$ \\
$\boxed{\Rightarrow   (\ell_j - u_j)Y_{ik} - \ell_kY_{ij} + \ell_i\ell_kx_j + u_j\ell_kx_i +\ell_i(u_j - \ell_j)x_k - \ell_iu_j\ell_k \leq 0 \quad  \toto \label{c6b}}$  \\
\noindent or symmetrically\\
$\ell_ix_jx_k + u_jx_ix_k + \ell_kx_ix_j - \ell_i\ell_kx_j - u_j\ell_kx_i - \ell_iu_jx_k + \ell_iu_j\ell_k \geq x_ix_jx_k \geq x_i(u_kx_j + u_jx_k - u_ju_k)$ \\
$- \ell_ix_jx_k +(u_k - \ell_k)x_ix_j + \ell_i\ell_kx_j +u_j(\ell_k - u_k)x_i + \ell_iu_jx_k - \ell_iu_j\ell_k \leq 0$ \\
$\boxed{\Rightarrow  - \ell_iY_{jk} +(u_k - \ell_k)Y_{ij} + \ell_i\ell_kx_j +u_j(\ell_k - u_k)x_i + \ell_iu_jx_k - \ell_iu_j\ell_k \leq 0 \quad  \toto \label{c6c}}$\\
\noindent or symmetrically\\
$\ell_ix_jx_k + u_jx_ix_k + \ell_kx_ix_j - \ell_i\ell_kx_j - u_j\ell_kx_i - \ell_iu_jx_k + \ell_iu_j\ell_k \geq x_ix_jx_k \geq x_k(u_ix_j + u_jx_i - u_iu_j)$ \\
$\Rightarrow (u_i - \ell_i) x_jx_k - \ell_kx_ix_j + \ell_i\ell_kx_j + u_j\ell_kx_i + u_j(\ell_i -u_i)x_k - \ell_iu_j\ell_k \leq 0$ \\
$\boxed{\Rightarrow (u_i - \ell_i) Y_{jk} - \ell_kY_{ij} + \ell_i\ell_kx_j + u_j\ell_kx_i +u_j(\ell_i -u_i)x_k - \ell_iu_j\ell_k \leq 0  \quad  \toto \label{c6d}}$\\
\noindent or
 $\ell_ix_jx_k + u_jx_ix_k + \ell_kx_ix_j - \ell_i\ell_kx_j - u_j\ell_kx_i - \ell_iu_jx_k + \ell_iu_j\ell_k \geq x_ix_jx_k \geq x_j(\ell_kx_i + \ell_ix_k - \ell_i\ell_k)$ \\
 $\Rightarrow - u_jx_ix_k  + u_j\ell_kx_i + \ell_iu_jx_k - \ell_iu_j\ell_k \leq 0 $\\
$\boxed{\Rightarrow   - u_jY_{ik}  + u_j\ell_kx_i + \ell_iu_jx_k - \ell_iu_j\ell_k \leq 0 \quad  \toto \label{c6e}}$\\
\noindent or
 $\ell_ix_jx_k + u_jx_ix_k + \ell_kx_ix_j - \ell_i\ell_kx_j - u_j\ell_kx_i - \ell_iu_jx_k + \ell_iu_j\ell_k \geq x_ix_jx_k \geq x_j(u_ix_k + u_kx_i - u_iu_k)$ \\
 $\boxed{\Rightarrow (u_i - \ell_i)Y_{jk} - u_jY_{ik} + (u_k  -\ell_k) Y_{ij} + (\ell_i\ell_k - u_iu_k)x_j + u_j\ell_kx_i + \ell_iu_jx_k - \ell_iu_j\ell_k \leq 0  \quad  \toto \label{c6f}}$\\

\begin{proposition}
    \label{prop6}
(T\ref{c6f}) cuts feasible solutions of $\QPR$, while~(T\ref{c6a})-(T\ref{c6e}) are redundant.
\end{proposition}
 
\begin{preuve}

\begin{enumerate}[(i)]
\item Consider the following solution  $(x,Y)$ satisfaying $h^t_{ij}(x,Y) \leq 0$, for all $(i,j,t) \in T$:
\begin{itemize}
\item $x_i = \frac{u_i + \ell_i}{2}$, $x_j = \frac{u_j + \ell_j}{2}$, and $x_k = \frac{u_k + \ell_k}{2}$
\item $Y_{ij} = \frac{u_iu_j + \ell_i\ell_j}{2}$, $Y_{ik} = \frac{\ell_i\ell_k + u_iu_k}{2}$, and $Y_{jk} = \frac{u_j\ell_k + \ell_ju_k}{2}$
\end{itemize}
Using this solution in Constraints~(T\ref{c6f}), we get: $\frac{u_i\ell_j\ell_k}{2} - \frac{\ell_iu_ju_k}{2} + \frac{u_iu_ju_k}{2}  + \frac{\ell_i\ell_ju_k}{2} - \frac{u_iu_j\ell_k}{2} -  \frac{\ell_i\ell_j\ell_k}{2} +  \frac{u_i\ell_ju_k}{2} + \frac{\ell_iu_j\ell_k}{2} = \frac{1}{2}(u_i - \ell_i)(\ell_j - u_j)(\ell_k - u_k) \geq 0$. This solution is thus cut off by the inequality~(T\ref{c6f}). 
\item (T\ref{c6a}) is equivalent to $(\ell_j - u_j)(Y_{ik} - \ell_kx_i - \ell_ix_k + \ell_i\ell_k) +  \ell_i(- Y_{jk}  + \ell_kx_j + \ell_jx_k -\ell_j\ell_k) \leq 0$ that is redundant with~(\ref{mc2}) and~(\ref{mc4}). By symmetry, the proof is similar for~(T\ref{c6b})-(T\ref{c6d}). 
\item (T\ref{c6e}) is equivalent to $u_j (- Y_{ik} + \ell_kx_i + \ell_ix_k - \ell_i\ell_k) \leq 0 $ that is redundant with~(\ref{mc4}). 
\end{enumerate}
\end{preuve}

\noindent \textbf{Family 7} We consider $(x_i -\ell_i)(x_j - \ell_j)(u_k - x_k) \geq 0$ and we get:\\
\noindent  $u_kx_ix_j + \ell_ix_jx_k + \ell_jx_ix_k - \ell_iu_kx_j - \ell_ju_kx_i - \ell_i\ell_jx_k + \ell_i\ell_ju_k \geq x_ix_jx_k \geq x_i(\ell_jx_k + \ell_kx_j - \ell_j\ell_k)$ \\
$\Rightarrow (\ell_k - u_k)x_ix_j -  \ell_ix_jx_k + \ell_iu_kx_j +\ell_j(u_k - \ell_k)x_i + \ell_i\ell_jx_k - \ell_i\ell_ju_k \leq 0$ \\
$\boxed{\Rightarrow (\ell_k - u_k)Y_{ij} -  \ell_iY_{jk} + \ell_iu_kx_j +\ell_j(u_k - \ell_k)x_i + \ell_i\ell_jx_k - \ell_i\ell_ju_k \leq 0\quad \toto \label{c7a}}$ \\ 
\noindent or symmetrically\\
$u_kx_ix_j + \ell_ix_jx_k + \ell_jx_ix_k - \ell_iu_kx_j - \ell_ju_kx_i - \ell_i\ell_jx_k + \ell_i\ell_ju_k \geq x_ix_jx_k \geq x_j(\ell_ix_k + \ell_kx_i - \ell_i\ell_k)$ \\
$\Rightarrow (\ell_k -u_k) x_ix_j -  \ell_jx_ix_k +\ell_i(u_k - \ell_k) x_j + \ell_ju_kx_i + \ell_i\ell_jx_k - \ell_i\ell_ju_k \leq 0$ \\
$\boxed{\Rightarrow (\ell_k -u_k) Y_{ij} -  \ell_jY_{ik} +\ell_i(u_k - \ell_k) x_j + \ell_ju_kx_i + \ell_i\ell_jx_k - \ell_i\ell_ju_k  \leq 0 \quad \toto \label{c7b}}$ \\ 
\noindent or \\
$u_kx_ix_j + \ell_ix_jx_k + \ell_jx_ix_k - \ell_iu_kx_j - \ell_ju_kx_i - \ell_i\ell_jx_k + \ell_i\ell_ju_k \geq x_ix_jx_k \geq x_i(u_jx_k + u_kx_j - u_ju_k)$ \\
$\Rightarrow - \ell_ix_jx_k + (u_j - \ell_j)x_ix_k + \ell_iu_kx_j + u_k(\ell_j -u_j)x_i + \ell_i\ell_jx_k - \ell_i\ell_ju_k \leq 0$ \\
$\boxed{\Rightarrow  - \ell_iY_{jk} + (u_j - \ell_j)Y_{ik} + \ell_iu_kx_j + u_k(\ell_j -u_j)x_i + \ell_i\ell_jx_k - \ell_i\ell_ju_k \leq 0 \quad \toto \label{c7c}}$ \\ 
\noindent or symmetrically\\
$u_kx_ix_j + \ell_ix_jx_k + \ell_jx_ix_k - \ell_iu_kx_j - \ell_ju_kx_i - \ell_i\ell_jx_k + \ell_i\ell_ju_k \geq x_ix_jx_k \geq x_j(u_ix_k + u_kx_i - u_iu_k)$ \\
$\Rightarrow (u_i - \ell_i)x_jx_k - \ell_jx_ix_k + u_k(\ell_i - u_i)x_j + \ell_ju_kx_i + \ell_i\ell_jx_k - \ell_i\ell_ju_k \leq 0$ \\
$\boxed{\Rightarrow   (u_i - \ell_i)Y_{jk} - \ell_jY_{ik} + u_k(\ell_i - u_i)x_j + \ell_ju_kx_i + \ell_i\ell_jx_k - \ell_i\ell_ju_k \leq 0 \quad \toto \label{c7d}}$ \\ 
\noindent or \\
$u_kx_ix_j + \ell_ix_jx_k + \ell_jx_ix_k - \ell_iu_kx_j - \ell_ju_kx_i - \ell_i\ell_jx_k + \ell_i\ell_ju_k \geq x_ix_jx_k \geq x_k(\ell_jx_i + \ell_ix_j - \ell_i\ell_j)$ \\
$\Rightarrow - u_kx_ix_j + \ell_iu_kx_j + \ell_ju_kx_i - \ell_i\ell_ju_k \leq 0$ \\
$\boxed{\Rightarrow   - u_kY_{ij} + \ell_iu_kx_j + \ell_ju_kx_i - \ell_i\ell_ju_k \leq 0 \quad \toto \label{c7e}}$ \\ 
\noindent or \\
$u_kx_ix_j + \ell_ix_jx_k + \ell_jx_ix_k - \ell_iu_kx_j - \ell_ju_kx_i - \ell_i\ell_jx_k + \ell_i\ell_ju_k \geq x_ix_jx_k \geq x_k(u_ix_j + u_jx_i - u_iu_j)$ \\
$\boxed{\Rightarrow - u_kY_{ij} + (u_i - \ell_i)Y_{jk} + (u_j - \ell_j)Y_{ik} + \ell_iu_kx_j + \ell_ju_kx_i + (\ell_i\ell_j - u_iu_j)x_k - \ell_i\ell_ju_k   \leq 0 \quad \toto \label{c7f}}$ \\

\begin{proposition}
    \label{prop7}
(T\ref{c7f}) cuts feasible solutions of $\QPR$, while~(T\ref{c7a})-(T\ref{c7e}) are redundant.
\end{proposition}
 
\begin{preuve}

\begin{enumerate}[(i)]
\item Consider the following solution  $(x,Y)$ satisfaying $h^t_{ij}(x,Y) \leq 0$, for all $(i,j,t) \in T$:
\begin{itemize}
\item $x_i = \frac{u_i + \ell_i}{2}$, $x_j = \frac{u_j + \ell_j}{2}$, and $x_k = \frac{u_k + \ell_k}{2}$
\item $Y_{ij} = \frac{\ell_iu_j + u_i\ell_j}{2}$, $Y_{ik} = \frac{\ell_i\ell_k + u_iu_k}{2}$, and $Y_{jk} = \frac{u_ju_k + \ell_j\ell_k}{2}$
\end{itemize}
Using this solution in Constraints~(T\ref{c7f}), we get $\frac{u_i\ell_j\ell_k}{2} - \frac{\ell_iu_ju_k}{2} + \frac{u_iu_ju_k}{2}  + \frac{\ell_iu_j\ell_k}{2} - \frac{u_i\ell_ju_k}{2} -  \frac{\ell_i\ell_j\ell_k}{2} +  \frac{u_iu_j\ell_k}{2} + \frac{\ell_i\ell_ju_k}{2} = \frac{1}{2}(u_i - \ell_i)(\ell_j - u_j)(\ell_k - u_k) \geq 0$. This solution is thus cut off by the inequality~(T\ref{c7f}). 
\item (T\ref{c7a}) is equivalent to $(\ell_k - u_k)(Y_{ij} - \ell_jx_i - \ell_ix_j + \ell_i\ell_j) +  \ell_i(- Y_{jk}  + \ell_jx_k + \ell_kx_j -\ell_j\ell_k) \leq 0$ that is redundant with~(\ref{mc2}) and~(\ref{mc4}). By symmetry, the proof is similar for~(T\ref{c7b})-(T\ref{c7d}).
\item (T\ref{c7e}) is equivalent to $u_k (- Y_{ij} + \ell_iix_j + \ell_jjx_i - \ell_i\ell_j) \leq 0 $ that is redundant with~(\ref{mc4}). 
\end{enumerate}
\end{preuve}

\noindent \textbf{Family 8} We consider $(x_i -\ell_i)(x_j - \ell_j)(x_k - \ell_k) \geq 0$ and we get:\\
\noindent $\ell_ix_jx_k + \ell_jx_ix_k + \ell_kx_ix_j - \ell_i\ell_kx_j - \ell_j\ell_kx_i - \ell_i\ell_jx_k + \ell_i\ell_j\ell_k \leq x_ix_jx_k \leq x_i(u_jx_k + \ell_kx_j - u_j\ell_k)$ \\
$\Rightarrow  \ell_ix_jx_k + (\ell_j - u_j)x_ix_k  - \ell_i\ell_kx_j +\ell_k(u_j - \ell_j)x_i - \ell_i\ell_jx_k + \ell_i\ell_j\ell_k \leq 0 $ \\
$\boxed{\Rightarrow \ell_iY_{jk} + (\ell_j - u_j)Y_{ik} - \ell_i\ell_kx_j +\ell_k(u_j - \ell_j)x_i - \ell_i\ell_jx_k + \ell_i\ell_j\ell_k \leq 0 \quad \toto \label{c8a}}$  \\
\noindent or symmetrically\\
$\ell_ix_jx_k + \ell_jx_ix_k + \ell_kx_ix_j - \ell_i\ell_kx_j - \ell_j\ell_kx_i - \ell_i\ell_jx_k + \ell_i\ell_j\ell_k \leq x_ix_jx_k \leq x_j(u_ix_k + \ell_kx_i - u_i\ell_k)$\\
$\Rightarrow (\ell_i -u_i)x_jx_k + \ell_jx_ix_k +\ell_k(u_i - \ell_i)x_j - \ell_j\ell_kx_i - \ell_i\ell_jx_k + \ell_i\ell_j\ell_k \leq 0$\\
 $\boxed{\Rightarrow (\ell_i -u_i)Y_{jk}+ \ell_jY_{ik} +\ell_k(u_i - \ell_i)x_j - \ell_j\ell_kx_i - \ell_i\ell_jx_k + \ell_i\ell_j\ell_k \leq 0 \quad \toto \label{c8b}}$  \\
\noindent or symmetrically\\
$\ell_ix_jx_k + \ell_jx_ix_k + \ell_kx_ix_j - \ell_i\ell_kx_j - \ell_j\ell_kx_i - \ell_i\ell_jx_k + \ell_i\ell_j\ell_k \leq x_ix_jx_k \leq x_k(u_jx_i + \ell_ix_j - \ell_iu_j)$\\
$\Rightarrow (\ell_j-u_j)x_ix_k + \ell_kx_ix_j - \ell_i\ell_kx_j - \ell_j\ell_kx_i +\ell_i(u_j - \ell_j)x_k + \ell_i\ell_j\ell_k \leq 0$ \\
$\boxed{\Rightarrow (\ell_j-u_j)Y_{ik}+ \ell_kY_{ij} - \ell_i\ell_kx_j - \ell_j\ell_kx_i +\ell_i(u_j - \ell_j)x_k + \ell_i\ell_j\ell_k \leq 0 \quad \toto \label{c8c} }$  \\
\noindent or  symmetrically\\
$\ell_ix_jx_k + \ell_jx_ix_k + \ell_kx_ix_j - \ell_i\ell_kx_j - \ell_j\ell_kx_i - \ell_i\ell_jx_k + \ell_i\ell_j\ell_k \leq x_ix_jx_k \leq x_i(\ell_jx_k + u_kx_j - \ell_ju_k)$ \\
$\Rightarrow \ell_ix_jx_k + (\ell_k -u_k)x_ix_j - \ell_i\ell_kx_j +\ell_j(u_k - \ell_k)x_i - \ell_i\ell_jx_k + \ell_i\ell_j\ell_k \leq 0$ \\
$\boxed{\Rightarrow \ell_iY_{jk} + (\ell_k -u_k)Y_{ij} - \ell_i\ell_kx_j +\ell_j(u_k - \ell_k)x_i - \ell_i\ell_jx_k + \ell_i\ell_j\ell_k \leq 0 \quad \toto \label{c8d} }$  \\
\noindent or symmetrically\\
$\ell_ix_jx_k + \ell_jx_ix_k + \ell_kx_ix_j - \ell_i\ell_kx_j - \ell_j\ell_kx_i - \ell_i\ell_jx_k + \ell_i\ell_j\ell_k \leq x_ix_jx_k \leq x_j(\ell_ix_k + u_kx_i - \ell_iu_k)$\\
$\Rightarrow  \ell_jx_ix_k + (\ell_k - u_k)x_ix_j + \ell_i(u_k - \ell_k)x_j - \ell_j\ell_kx_i - \ell_i\ell_jx_k + \ell_i\ell_j\ell_k \leq 0 $\\
$\boxed{\Rightarrow  \ell_jY_{ik} + (\ell_k - u_k)Y_{ij} + \ell_i(u_k - \ell_k)x_j - \ell_j\ell_kx_i - \ell_i\ell_jx_k + \ell_i\ell_j\ell_k \leq 0 \quad \toto \label{c8e}}$ \\
\noindent or symmetrically\\
$\ell_ix_jx_k + \ell_jx_ix_k + \ell_kx_ix_j - \ell_i\ell_kx_j - \ell_j\ell_kx_i - \ell_i\ell_jx_k + \ell_i\ell_j\ell_k \leq x_ix_jx_k \leq x_k(\ell_jx_i + u_ix_j - u_i\ell_j)$\\
$\Rightarrow (\ell_i -u_i)x_jx_k + \ell_kx_ix_j - \ell_i\ell_kx_j - \ell_j\ell_kx_i +\ell_j(u_i - \ell_i)x_k + \ell_i\ell_j\ell_k \leq  0 $\\
$\boxed{\Rightarrow  (\ell_i -u_i)Y_{jk} + \ell_kY_{ij} - \ell_i\ell_kx_j - \ell_j\ell_kx_i +\ell_j(u_i - \ell_i)x_k + \ell_i\ell_j\ell_k  \leq 0 \quad \toto \label{c8f}}$ \\

\begin{proposition}
    \label{prop8}
(T\ref{c8a})-(T\ref{c8f}) are redundant for $\QPR$.
\end{proposition}
 
\begin{preuve}
(T\ref{c8a}) is equivalent to $(\ell_j - u_j)(Y_{ik} - \ell_kx_i - \ell_ix_k + \ell_i\ell_k) +  \ell_i(Y_{jk}  - \ell_kx_j - u_jx_k + u_j\ell_k) \leq 0$ that is redundant with~(\ref{mc1}) and~(\ref{mc4}). By symmetry, the proof is similar for~(T\ref{c8b})-(T\ref{c8f}).
\end{preuve}

\bigskip

\setcounter{equation}{7}

We now  define the set $\setIndT' = \{(i,j,k,t): (i,j,k) \in \setT, t=1, \ldots, 12 \}$, and from Propositions~\ref{prop1}--\ref{prop8}, we introduce the set of general triangle inequalities,
$$\setGT = \Big \{(x,Y) \, : \, h^{t'}_{ijk}(x,Y) \leq 0\quad  \forall (i,j,k,t) \in \setIndT'  \Big \} $$
with  $h^{t'}_{ijk}(x,Y)$:
\begin{small}
\begin{numcases}{} 
(\ell_k -u_k)Y_{ij} + (\ell_j - u_j)Y_{ik} - u_i Y_{jk}+ u_iu_kx_j + (u_ju_k - \ell_j\ell_k)x_i + u_iu_jx_k - u_iu_ju_k  & $t=1$ \label{t1}  \\ 
(\ell_k -u_k)Y_{ij} - u_jY_{ik}+ ( \ell_i -u_i)Y_{jk} + (u_iu_k - \ell_i\ell_k)x_j + u_ju_kx_i + u_iu_jx_k - u_iu_ju_k  &$t=2$\label{t2}  \\ 
 - u_kY_{ij} + (\ell_j - u_j)Y_{ik}+ ( \ell_i -u_i)Y_{jk} + u_iu_kx_j + u_ju_kx_i + ( u_iu_j - \ell_i\ell_j)x_k - u_iu_ju_k  &$t=3$ \label{t3}  \\ 
 (u_j - \ell_j) Y_{ik} + (\ell_k - u_k) Y_{ij} + u_iY_{jk}+ (\ell_ju_k -u_j\ell_k)x_i - u_i\ell_kx_j - u_iu_jx_k + u_iu_j\ell_k  &$t=4$\label{t4}  \\ 
 u_j Y_{ik} + (\ell_k - u_k) Y_{ij} + (u_i-\ell_i)Y_{jk} -u_j\ell_kx_i + (\ell_iu_k - u_i\ell_k)x_j - u_iu_jx_k + u_iu_j\ell_k  0 &$t=5$\label{t5}  \\ 
(u_k - \ell_k) Y_{ij} + u_iY_{jk} + (\ell_j - u_j) Y_{ik} -u_iu_kx_j +(u_j\ell_k - \ell_ju_k)x_i - u_i\ell_jx_k + u_i\ell_ju_k & $t=6$\label{t6}  \\ 
u_kY_{ij} + (u_i - \ell_i) Y_{jk} + (\ell_j - u_j) Y_{ik} -u_iu_kx_j - \ell_ju_kx_i +(\ell_iu_j - u_i\ell_j)x_k + u_i\ell_ju_k & $t=7$\label{t7}  \\ 
(u_k - \ell_k)Y_{ij} + (\ell_i - u_i) Y_{jk} + u_jY_{ik} +(u_i\ell_k - \ell_iu_k)x_j - u_ju_kx_i - \ell_iu_jx_k + \ell_iu_ju_k \leq & $t=8$\label{t8}  \\ 
u_kY_{ij} + (\ell_i - u_i) Y_{jk} + (u_j - \ell_j)Y_{ik} - \ell_iu_kx_j - u_ju_kx_i + (u_i\ell_j - \ell_iu_j) x_k + \ell_iu_ju_k   & $t=9$\label{t9}  \\ 
- u_i Y_{jk} + (u_j - \ell_j)Y_{ik} + (u_k - \ell_k) Y_{ij} + u_i\ell_kx_j + (\ell_j\ell_k - u_ju_k)x_i + u_i\ell_jx_k - u_i\ell_j\ell_k & $t=10$\label{t10}  \\ 
 (u_i - \ell_i)Y_{jk} - u_jY_{ik} + (u_k  -\ell_k) Y_{ij} + (\ell_i\ell_k - u_iu_k)x_j + u_j\ell_kx_i + \ell_iu_jx_k - \ell_iu_j\ell_k  & $t=11$\label{t11}  \\ 
- u_kY_{ij} + (u_i - \ell_i)Y_{jk} + (u_j - \ell_j)Y_{ik} + \ell_iu_kx_j + \ell_ju_kx_i + (\ell_i\ell_j - u_iu_j)x_k - \ell_i\ell_ju_k    & $t=12$\label{t12} 
\end{numcases}
\end{small}

We now state in Proposition~\ref{propTriang01} that these inequalities amounts to the Triangle inequalities, introduced in~\cite{Pad89}, when for all $i \in \setI$, $\ell_i = 0$ and $u_i =1$. 

\begin{proposition}
  \label{propTriang01}
For all $(i,j,k,t) \in \setIndT'$, inequalities $h^{t'}_{ijk}(x,Y) \leq 0$ are an extension of the Triangle inequalities to the  case of general upper and lower bounds (i.e. $x_i \in [\ell_i,u_i]$). 
\end{proposition}
\begin{preuve}
By setting in inequalities $h^{t'}_{ijk}(x,Y) \leq 0$ each lower and upper bound by values $0$ and $1$ respectively, we come back to the classical Triangle inequalities.
\end{preuve}

\section{Computing a strenghtened quadratic convex relaxation}
\label{sec:bestrelax}

By adding the general triangle inequalities to $\QPR$, we obtain a strenghened family of quadratic convex relaxation to $\QP$:
\begin{numcases}{\QPT}
 \min \;\;\;  \langle S_0, xx^T \rangle +c_0^Tx + \langle Q_0 - S_0 , Y \rangle \nonumber \\
 \mbox{s.t.} \,\,\,\,  (\ref{bound})(\ref{real})  \nonumber \\
 \qquad  \langle S_r, xx^T \rangle +c_r^Tx + \langle Q_r - S_r , Y \rangle  \leq b_r  & $r=1, \ldots, m$\,\, \nonumber \\
 \qquad h^t_{ij}(x,Y) \leq 0 & $(i,j,t) \in \setIndT $\nonumber \\
\qquad h^{t'}_{ij}(x,Y) \leq 0 & $(i,j,k,t) \in \setIndT' $\nonumber \\
 \qquad Y_{jj} = Y_{ij} & $ (i,j) \in \overline{\setIMC}$ \nonumber
\end{numcases}

\medskip

\par We then consider the problem of finding the best set of positive semi-definite matrices $S_0, \ldots ,S_m$, in the sense that the optimal solution value of $\QPT$ is as large as possible. This amounts to solving the following problem $\OPT$:
\begin{numcases}{\OPT}
 \max_{\begin{array}{c} \scriptstyle{S_0, \ldots, S_m \succeq 0}\end{array}}
  v\QPT  \nonumber
\end{numcases}

\par where $v\QPT$ is the optimal value of problem $\QPT$. We further prove that $v\OPT$ is equal to the optimal value of the semi-definite relaxation of $\QP$ called "Shor's plus RLT plus Triangle".

\begin{theorem}
\label{theosdp}
Let $(SDP)$ be the  "Shor's plus RLT plus Triangle" semi-definite relaxation:
\begin{numcases}{(SDP)}
\min  f(X,x) \equiv \langle Q_0, X \rangle  +  c_0^T x  \nonumber \\
\mbox{s.t.} \,\,\,\, (X,x) \in \setSDP  \nonumber\\
\qquad h^t_{ij} (X,x) \equiv \langle M^t_{ij} , X \rangle +  (v^t_{ij})^Tx + l^t_{ij}  \leq 0, & $(i,j,t) \in \overline{\setIndT}$ \nonumber\\
\qquad h^{'t}_{ijk} (X,x) \equiv \langle M^{'t}_{ijk} , X \rangle +  (v^{'t}_{ijk})^Tx + l^{'t}_{ijk}\leq 0, & $(i,j,k,t) \in \setIndT'$ \nonumber
\end{numcases}
\noindent where  $\overline{\setIndT} = \{(i,j,t): (i,j) \in \overline{\setIMC}, t=1, \ldots, 4 \}$, and $\forall (i,j,t) \in \overline{\setIndT}$, matrices $M^t_{ij}$, vectors  $v^t_{ij}$, and scalars  $l^t_{ij}$ are the coefficients of constraints $h^t_{ij}(X,x) \leq 0$, $\forall (i,j,k,t) \in \setIndT'$, matrices $M^{'t}_{ijk}$, vectors  $v^{'t}_{ijk}$, and scalars  $l^{'t}_{ijk}$ are the coefficients of constraints $h^{'t}_{ijk}(X,x) \leq 0$, and $\setSDP$ is the following set:
\begin{numcases}{\setSDP=(x,X)}
 \langle Q_{r}, X\rangle  + c^T_{r}x  \leq b_r &  $r \in \setR$  \label{sdp_q1}  \\
X_{ii} -  u_ix_i -\ell_ix_i + u_j\ell_i \leq 0  & $i \in \setI$ \label{sdp_mcii1} \\ 
- X_{ii} + 2u_i x_i -u^2_i \leq 0 &$i \in \setI$ \label{sdp_mcii2}\\
- X_{ii}  + 2\ell_ix_i - \ell^2_i \leq 0 & $i \in \setI$\label{sdp_mcii3}\\
 \left ( \begin{array}{ll}
   1 & x^T \\
   x & X 
 \end{array}\right ) \succeq 0  \label{sdp_1}  \\
 x \in \mathbb{R}^n \quad X \in \mathcal{S}_n   \label{sdp_2} 
\end{numcases}
  It holds that $v\OPT = v(SDP)$.
\end{theorem}

\begin{preuve} 

\noindent $\diamond$ To prove that $v\OPT \leq v(SDP)$, we show that $v(\QPTbar) \leq v(SDP)$ for any  $\bar{S}_0 ,\ldots, \bar{S}_m \in \mathcal{S}_n^+$, which in turn implies that $v\OPT \leq v(SDP)$ since the right hand side is constant. For this, we show that, if $(\bar{x},\bar{X})$ is feasible for $(SDP)$, then $(x,Y):= (\bar{x},\bar{X})$ is $i)$ feasible for $(\QPTbar)$ and $ii)$ its objective value is less or equal than $v(SDP)$. Since $(\QPTbar)$ is a minimization problem, $v(\QPTbar)\leq v(SDP)$ follows.
 
\begin{enumerate}[i)]
\item   We prove that $(x,Y)$ is  feasible  to $(\QPTbar)$. Constraints of sets $\setMC$ and $\setGT$ are obviously satisfied. We now prove that Constraints~(\ref{cont}) are satisfied:
\begin{align*}
\langle  \bar{S}_r , xx^T  \rangle +c_r^Tx + \langle Q_r - \bar{S}_r , Y \rangle & =  \langle  \bar{S}_r , \bar{x}\bar{x}^T  \rangle +c_r^T\bar{x} + \langle Q_r - \bar{S}_r , \bar{X} \rangle \\
& = \langle \bar{S}_r , \bar{x}\bar{x}^T - \bar{X} \rangle + c_r^T\bar{x} + \langle Q_r  , \bar{X} \rangle \leq  b_r
\end{align*}
from Constraints~(\ref{sdp_q1}) and ~(\ref{sdp_1}), and since  $\bar{S}_r \succeq 0$.

\item  Let us compare the objective values. For this, we prove that   $\langle \bar{S}_0, \bar{x}\bar{x}^T\rangle  +c_0^T\bar{x} + \langle Q_0 -   \bar{S}_0 , \bar{X} \rangle - \langle Q_0 , \bar{X} \rangle - c_0^T\bar{x}  \leq 0$ or that $ \langle \bar{S}_0 , \bar{x}\bar{x}^T -\bar{X}  \rangle   \leq 0$. This last inequality follows from  $\bar{S}_0 \succeq 0$ and Constraint~(\ref{sdp_1}). \\
\end{enumerate}
\bigskip

\noindent  $\diamond$ Let us secondly prove that $v\OPT \geq v(SDP)$ or equivalently  $v\OPT \geq v(D)$ where $(D)$ is the dual of $(SDP)$:
\begin{small}
\begin{numcases}{(D)}
  \max \, g(\alpha,\Phi,\Delta,\rho) =  - \displaystyle{\sum_{r=1}^m \alpha_r} b_r + \displaystyle{\sum_{\tiny{(i,j) \in \setIMC}}} \displaystyle{\sum_{t=1}^{4}} \phi_{ij}^{t} l_{ij}^{t} +  \displaystyle{\sum_{\tiny{(i,j,k) \in \setT}}} \displaystyle{\sum_{t=1}^{12}} \delta_{ijk}^{t} l_{ijk}^{'t} \nonumber + u^T\ell\varphi^1 - u^Tu\varphi^2 - \ell^T\ell \varphi^3 - \rho\\
\mbox{s.t.} \nonumber \\
 \quad  S= Q_0 + \displaystyle{\sum_{r=1}^m} \alpha_r Q_r + \Phi + \Delta  \label{dsdp_1} \\ 
 \quad  d= c_0 + \displaystyle{\sum_{r=1}^m \alpha_r} c_r +   \displaystyle{\sum_{\tiny{(i,j) \in \setIMC}}} \displaystyle{\sum_{t=1}^{4}} \phi_{ij}^{t} v_{ij}^{t} + \displaystyle{\sum_{\tiny{(i,j,k) \in \setT}}} \displaystyle{\sum_{t=1}^{12}} \delta_{ijk}^{t} v_{ijk}^{'t}   - \varphi^{1T}(u + l) + 2\varphi^{2T}u + 2\varphi^{3T}\ell \label{dsdp_2} \\ 
 \quad \Phi = \displaystyle{\sum_{\tiny{(i,j) \in \setIMC}}} \displaystyle{\sum_{t=1}^{4}} \phi_{ij}^{t} M_{ij}^{t} + diag(\varphi^{1} - \varphi^{2}- \varphi^{3}) \label{dsdp_3}  \\
  \quad \Delta = \displaystyle{\sum_{\tiny{(i,j,k) \in \setT}}} \displaystyle{\sum_{t=1}^{12}} \delta_{ijk}^{t} M_{ijk}^{'t}  \label{dsdp_4} \\
\quad  \left ( \begin{array}{ll}
   \rho & \frac{1}{2}d^T \\
   \frac{1}{2}d & S
     \end{array}\right ) \succeq 0  \label{dsdp_5}  \\
  \quad \alpha \in \mathbb{R}^m_+, \,\, \varphi^1, \varphi^2, \varphi^3  \in \mathbb{R}^n_+, \,\, \Phi \in \mathcal{S}_n, \,\, \phi^t_{ij} \geq 0 \,  (i,j,t) \in \overline{T}, \,\,   \Delta \in \mathcal{S}_n, \,\, \delta^t_{ijk} \geq 0  \,  (i,j,k,t) \in T'
\end{numcases}
\end{small}
\noindent where $\alpha \in \mathbb{R}^m_+$ are the dual variables associated to constraints~(\ref{sdp_q1}), and $\phi^t_{ij}$ are the dual variable associated with Constraints~(\ref{mc1})--(\ref{mc4}), respectively, $\varphi^t$, $t=1,\ldots,3$ are the dual variables associated to Constraints~(\ref{sdp_mcii1})--(\ref{sdp_mcii3}), respectively, and  $\delta^t_{ijk}$ are the dual variables associated with Constraints~(\ref{t1})--(\ref{t12}), respectively.\\

\par Let $(\alpha^*,\Phi^*,\Delta^*,\rho^*)$ be an optimal solution to $(D)$, we build the following positive semi-definite matrices:  $\forall r \in\setR$, $\bar{S}_r = \mathbf{0}_n$, and $\bar{S}_0  = S^* =  Q_0 + \displaystyle{\sum_{r=1}^m} \alpha^*_r Q_r  + \Phi^* + \Delta^*$. By Constraint~(\ref{dsdp_5}), $\bar{S}_0 \succeq 0$, and $(\bar{S}_0, \ldots ,\bar{S}_m)$ forms a feasible solution to $\OPT$. The objective value of this solution is equal to $v(\QPTbar)$.

We now prove that $v(\QPTbar) \geq v(D)$. For this, we prove that for any feasible solution $(\bar{x},\bar{Y})$ to $(\QPTbar)$, the associated objective value is not smaller than $g(\alpha^*,\Phi^*,\Delta^*,\rho^*)$. Denote by  $\theta$ the difference between the objective values, i.e., $\theta =  \langle  \bar{S}_0 ,\bar{x}\bar{x}^T  \rangle +c_0^T\bar{x} + \langle Q_0 - \bar{S}_0 , \bar{Y} \rangle - g(\alpha^*,\Phi^*,\Delta^*,\rho^*)$. We below prove that $\theta \geq 0$.
\begin{small}
\begin{align*}
  \theta &= \langle  \bar{S}_0 , \bar{x}\bar{x}^T  \rangle +c_0^T\bar{x} + \langle Q_0 - \bar{S}_0 , \bar{Y} \rangle + \displaystyle{\sum_{r=1}^m \alpha^*_r} b_r - \displaystyle{\sum_{\tiny{(i,j) \in \setIMC}}} \displaystyle{\sum_{t=1}^{4}} \phi_{ij}^{t*} l_{ij}^{t} -  \displaystyle{\sum_{\tiny{(i,j,k) \in \setT}}} \displaystyle{\sum_{t=1}^{12}} \delta_{ijk}^{t*} l_{ijk}^{'t}  - u^T\ell\varphi^{1*} + u^Tu\varphi^{2*}\\
  &\quad + \ell^T\ell \varphi^{3*} + \rho^*\\
  \theta & = \langle  \bar{S}_0 , \bar{x}\bar{x}^T  \rangle+  c_0^T\bar{x} - \langle \displaystyle{\sum_{r=1}^m} \alpha^*_r Q_r  + \Phi^* + \Delta^* ,  \bar{Y} \rangle + \displaystyle{\sum_{r=1}^m \alpha^*_r} b_r - \displaystyle{\sum_{\tiny{(i,j) \in \setIMC}}} \displaystyle{\sum_{t=1}^{4}} \phi_{ij}^{t*} l_{ij}^{t} -  \displaystyle{\sum_{\tiny{(i,j,k) \in \setT}}} \displaystyle{\sum_{t=1}^{12}} \delta_{ijk}^{t*} l_{ijk}^{'t}  - u^T\ell\varphi^{1*}\\
  & \quad + u^Tu\varphi^{2*} + \ell^T\ell \varphi^{3*}  + \rho^*\\
& \mbox{  since }   Q_0 - \bar{S}_0 =  -(\displaystyle{\sum_{r=1}^m} \alpha^*_r Q_r  + \Phi^* + \Delta^*)\\
\theta  & =   \langle  \bar{S}_0 , \bar{x}\bar{x}^T  \rangle + c_0^T\bar{x} + \displaystyle{\sum_{r=1}^m} \alpha^*_r ( b_r - \langle Q_r ,  \bar{Y} \rangle) - \langle \Phi^* + \Delta^*, \bar{Y} \rangle  - \displaystyle{\sum_{\tiny{(i,j) \in \setIMC}}} \displaystyle{\sum_{t=1}^{4}} \phi_{ij}^{t*} l_{ij}^{t} -  \displaystyle{\sum_{\tiny{(i,j,k) \in \setT}}} \displaystyle{\sum_{t=1}^{12}} \delta_{ijk}^{t*} l_{ijk}^{'t} - u^T\ell\varphi^{1*} \\
  & \quad + u^Tu\varphi^{2*} + \ell^T\ell \varphi^{3*} + \rho^*\\
\theta  & \geq   \langle  \bar{S}_0 , \bar{x}\bar{x}^T  \rangle + c_0^T\bar{x} + \displaystyle{\sum_{r=1}^m} \alpha^*_r c^T_r \bar{x} - \langle \Phi^* + \Delta^*, \bar{Y} \rangle -\displaystyle{\sum_{\tiny{(i,j) \in \setIMC}}} \displaystyle{\sum_{t=1}^{4}} \phi_{ij}^{t*} l_{ij}^{t} -  \displaystyle{\sum_{\tiny{(i,j,k) \in \setT}}} \displaystyle{\sum_{t=1}^{12}} \delta_{ijk}^{t*} l_{ijk}^{'t}  - u^T\ell\varphi^{1*} + u^Tu\varphi^{2*} \\
  & \quad + \ell^T\ell \varphi^{3*} + \rho^*
&\intertext{as $ c_r^T \bar{x} + \langle Q_r ,\bar{Y} \rangle \leq  b_r $ and $\alpha^*_r \geq 0$. Moreover, by Constraint~(\ref{dsdp_3}) and~(\ref{dsdp_4}) we get:}
\theta   & \geq   \langle  \bar{S}_0 , \bar{x}\bar{x}^T  \rangle + c_0^T\bar{x} + \displaystyle{\sum_{r=1}^m} \alpha^*_r c^T_r \bar{x} -  \displaystyle{\sum_{\tiny{(i,j) \in \setIMC}}} \displaystyle{\sum_{t=1}^{4}}  \langle \phi_{ij}^{t*} M_{ij}^{t}, \bar{Y} \rangle  - \displaystyle{\sum_{\tiny{(i,j) \in \setIMC}}} \displaystyle{\sum_{t=1}^{4}} \phi_{ij}^{t*} l_{ij}^{t} -  \displaystyle{\sum_{\tiny{(i,j,k) \in \setT}}} \displaystyle{\sum_{t=1}^{12}} \langle \delta_{ijk}^{t*} M_{ijk}^{'t} , \bar{Y} \rangle \\
& \quad -  \displaystyle{\sum_{\tiny{(i,j,k) \in \setT}}} \displaystyle{\sum_{t=1}^{12}} \delta_{ijk}^{t*} l_{ijk}^{'t}  - \langle diag(\varphi^{1*} - \varphi^{2*}- \varphi^{3*}) , \bar{Y} \rangle - u^T\ell\varphi^{1*} + u^Tu\varphi^{2*} + \ell^T\ell \varphi^{3*} + \rho^* \\
\intertext{By Constraints~(\ref{mc1})--(\ref{mc4}) and~(\ref{t1})--(\ref{t12}), and since all coefficients $\bar{\phi}^t_{ij}$, $\bar{\delta}^t_{ijk}$ and $\bar{\varphi^t_i}$ are non-negative, we get:}
\theta  & \geq  \langle  \bar{S}_0 , \bar{x}\bar{x}^T  \rangle+   c^T_0\bar{x} + \displaystyle{\sum_{r=1}^m \alpha_r} c^T_r\bar{x}  - (\varphi^{1T*}(u + l) - 2\varphi^{2T*}u - 2\varphi^{3T*}\ell)^T\bar{x}  +   \displaystyle{\sum_{\tiny{(i,j) \in \setIMC}}} \displaystyle{\sum_{t=1}^{4}} \phi_{ij}^{t*} v_{ij}^{tT}\bar{x} \\
&\quad + \displaystyle{\sum_{\tiny{(i,j,k) \in \setT}}} \displaystyle{\sum_{t=1}^{12}} \delta_{ijk}^{t*} v_{ijk}^{'tT}\bar{x}  + \rho^* \\
\theta  & \geq   (c_0 + \displaystyle{\sum_{r=1}^m \alpha_r} c_r +   \displaystyle{\sum_{\tiny{(i,j) \in \setIMC}}} \displaystyle{\sum_{t=1}^{4}} \phi_{ij}^{t*} v_{ij}^{t} + \displaystyle{\sum_{\tiny{(i,j,k) \in \setT}}} \displaystyle{\sum_{t=1}^{12}} \delta_{ijk}^{t*} v_{ijk}^{'t} - \varphi^{1T*}(u + l) + 2\varphi^{2T*}u + 2\varphi^{3T*}\ell)^T\bar{x} \\
& \quad + \langle  \bar{S}_0 , \bar{x}\bar{x}^T  \rangle + \rho^*  \\
\theta & \geq \langle  \bar{S}_0 , \bar{x}\bar{x}^T  \rangle + d^{*T}\bar{x}+ \rho^* \textrm{ by Constraint~(\ref{dsdp_2})}
\end{align*}
\end{small}
We end the proof by showing that $\langle  \bar{S}_0 , \bar{x}\bar{x}^T  \rangle + d^{*T}\bar{x}+ \rho^* \geq 0$. From Constraint~(\ref{dsdp_5}), we know that for all $x \in \mathbb{R}^n$, \begin{small}$ \left ( \begin{array}{ll}    1 \\   x      \end{array}\right )^T \left ( \begin{array}{ll}    \rho^* & \frac{1}{2}d^{*T} \\   \frac{1}{2}d^* &  \bar{S}_0     \end{array}\right )   \left ( \begin{array}{ll}    1 \\   x      \end{array}\right ) \geq 0$ \end{small}, which prove that $\theta \geq 0$.
 \end{preuve}

\medskip

From the proof of Theorem~\ref{theosdp}, we can caracterize a set of optimal matrices  $(S_0^*, \ldots , S_m^*)$.
\begin{corollary}
  \label{corsdp}
 The following positive semi-definite matrices allow to build an optimal solution $(S_0^*, \ldots , S_m^*)$ of $\OPT$:
\begin{enumerate}[i)]
\item $\forall r = 1 , \ldots , m$, $S_r^* = \mathbf{0}_n$
\item $S_0^* =  Q_0 + \displaystyle{\sum_{r=1}^m} \alpha^*_r Q_r + \Phi^* + \Delta^* $,  where:\\
  \noindent $\diamond$ $\alpha^*$ is the vector of optimal dual variables associated with Constraints~(\ref{sdp_q1}),\\
  \noindent $\diamond$  matrix $\Phi^* = \displaystyle{\sum_{\tiny{(i,j) \in \setIMC}}} \displaystyle{\sum_{t=1}^{4}} \phi_{ij}^{t*} M_{ij}^{t*} + diag(\varphi^{1*} - \varphi^{2*}- \varphi^{3*})$, where  $\varphi^{1*}$, $\varphi^{2*}$, $\varphi^{3*}$ are the vectors of optimal dual variables associated with Constraints~(\ref{sdp_mcii1})--(\ref{sdp_mcii3}),  and $\forall (i,j,t) \in \overline{\setIndT}$, $\phi^{t*}_{ij}$ is the optimal dual variables associated with Constraints~(\ref{mc1})--(\ref{mc4}), respectively.\\
 \noindent $\diamond$ matrix  $\Delta^* = \displaystyle{\sum_{\tiny{(i,j,k) \in \setT}}} \displaystyle{\sum_{t=1}^{12}} \delta_{ijk}^{t*} M_{ijk}^{'t*}$, where $\forall (i,j,k,t) \in \setIndT'$, $\delta_{ijk}^{t*}$ is the optimal dual variables associated with Constraints~(\ref{t1})--(\ref{t12}), respectively.
\end{enumerate}
\end{corollary}
\medskip

To sum up, we obtain the following quadratic convex relaxation to $\QP$:
\begin{numcases}{(P^*)}
 \min \;\;\; f_{0,S_0^*}(x,Y) = \langle Q_0 + \displaystyle{\sum_{r=1}^m} \alpha^*_r Q_r + \Phi^* + \Delta^*, xx^T \rangle +c_0^Tx + \langle - \displaystyle{\sum_{r=1}^m} \alpha^*_r Q_r - \Phi^* - \Delta^* , Y \rangle \nonumber \\
 \mbox{s.t.} \,\,\,\,   f_r(x,Y) = \langle Q_r, Y \rangle +c_r^Tx  \leq b_r &$r \in \setR$  \nonumber\\
\qquad h^t_{ij}(x,Y) \leq 0 & $(i,j,t) \in \setIndT $\nonumber \\
\qquad h^{t'}_{ij}(x,Y) \leq 0 & $(i,j,k,t) \in \setIndT' $\nonumber \\
 \qquad Y_{jj} = Y_{ij} & $ (i,j) \in \overline{\setIMC}$ \nonumber
\end{numcases}
\medskip

As stated in Theorem~\ref{theosdp}, the optimal value of $(P^*)$ is equal to the optimal value of $(SDP)$, and we now by Proposition~\ref{prop1}--\ref{prop8} that this relaxation is tighter than the "shor plus RLT" relaxation. This sharp relaxation can then be used within a \bb~algorithm to solve $\QP$ to global optimality.

\section{Using the dynamic bundle method for heuristically solving $(SDP)$ and for separating inequalities}
\label{sec:bundle}

In this section we propose to separate the set of  inequalities $\setMC \cup \setGT$ by heuristically solving $(SDP)$, thanks to a dynamic bundle method. For this, we introduce a parameter $p$ that controls the number of considered constraints in $(SDP)$, and thus in the associated computed quadratic convex relaxation. Following the idea of~\cite{BELW17}, we design sub-gradient algorithm within a Lagrangian duality framework using the callable Conic Bundle library of~\cite{cb}.

\medskip
\par To describe the algorithm, we consider $(SDP)$ as a maximization problem, by changing the sign of its objective function. We then consider a partial Lagrangian dual of  $(SDP)$ where we dualize the set of constraints $\setMC \cup \setGT$.  For this, with each constraint $h^t_{ij}(x,X) \leq 0$, we associate a non-negative Lagrange multiplier $\phi^t_{ij}$, and with each constraints $h^{'t}_{ijk}(x,X) \leq 0$ a non-negative Lagrange multiplier $\delta^t_{ijk}$. We now consider the partial Lagrangian: 
$$\mathcal{L}_{\overline{\setIndT},\setIndT'}(x,X,\phi,\delta) = - \langle Q_0, X \rangle  -  c_0^T x - \displaystyle{\sum_{(i,j,t) \in \overline{\setIndT}}} \phi^t_{ij} h^t_{ij}(x,X)- \displaystyle{\sum_{(i,j,k,t) \in \setIndT'}} \delta^t_{ijk} h^{'t}_{ijk}(x,X)$$
and we obtain the dual functional
$$g_{\overline{\setIndT},\setIndT'}(\phi,\delta) = \displaystyle{\max_{\scriptsize{(x,X) \in \setSDP}}} \quad \mathcal{L}_{\overline{\setIndT},\setIndT'}(x,X,\phi,\delta)$$
By minimizing this dual functional we obtain the partial Lagrangian dual problem:
\begin{numcases}{(LD_{\overline{\setIndT},\setIndT'})}
  \min_{\begin{tiny}\begin{array}{c}\phi^t_{ij} \geq 0,\; (i,j,t) \in \overline{\setIndT}  \\ \delta^t_{ijk} \geq 0,\; (i,j,k,t) \in \setIndT'   \end{array}\end{tiny}}  g_{\overline{\setIndT},\setIndT'}(\phi,\delta)  \nonumber 
\end{numcases}

Problem $(LD_{\overline{\setIndT},\setIndT'})$ can then be solved with the bundle method for which a detailed description is available in~\cite{fischeretal}.  However, the number of elements in $\overline{T} \cup T'$ is $4{n \choose 2} + 12{n \choose 3}$, and we are interested only in the subset of $\overline{T} \cup T'$ for which the constraints $h_{ij}^t(x,X) \leq 0$ and $h_{ijk}^{'t}(x,X) \leq 0$ are active at the optimum. In order to preserve efficiency we adopt another idea from~\cite{fischeretal} that consists in dynamically adding and removing constraints in the course of the algorithm. Then, we now consider $\overline{\mathcal{T}} \cup \mathcal{T}' \subseteq \overline{T} \cup T'$ and work with the function:
$$g_{\overline{\mathcal{T}},\mathcal{T'}}(\phi,\delta) = \displaystyle{\max_{\scriptsize{(x,X) \in \setSDP}}} \quad \mathcal{L}_{\overline{\mathcal{T}},\mathcal{T'}}(x,X,\phi,\delta).$$
Initially we set $\overline{\mathcal{T}}\cup \mathcal{T'}=\emptyset$ and after a first function evaluation we separate violated inequalities and add the elements to set $\overline{\mathcal{T}} \cup \mathcal{T'}$ accordingly. We keep on updating this set in course of the bundle iterations by removing elements with associated multiplier close to zero and separate newly violated constraints. Convergence for dynamic bundle methods has been analyzed in detail in~\cite{BeSa:09}, giving a positive answer for convergence properties in a rather general setting.
\bigskip

\par In our context, we know that any feasible dual solution to $(SDP)$ allows us to build a quadratic convex relaxation of $\QP$. Better this solution is, sharper is the associated bound at the root node of the \bb~process. Another idea to reduce more the solution time of $(SDP)$ is to consider a parameter $p$ that is an upper bound on the cardinality of $\overline{\mathcal{T}} \cup \mathcal{T'}$ (i.e. $\lvert \overline{\mathcal{T}} \cup \mathcal{T'} \rvert \leq p$). In other words, $p$ is the maximum number of constraints considered in the reduced problem. Introducing this parameter $p$ leads to a dual heuristic that has two extreme cases: 
\begin{itemize}
\item if $p= 4{n \choose 2} + 12{n \choose 3} $, we solve $(SDP)$ and get the associated dual solution as described in Corollary~\ref{corsdp}.
\item if $p=0$,  we make a single iteration: we get the optimal solution of the reduced problem obtained from $(SDP)$ where we drop all constraints~of sets $\setMC$ and $\setGT$, this amounts to solving the "Shor's plus diagonal RLT" semi-definite relaxation.
\end{itemize}

\section{Computational results}
\label{sec:exp}
In this section, we compare our algorithm \miqcrt~with \glomiqo~(\cite{MSF15}),  \baron~(\cite{baron}), and the original \miqcr~method~(\cite{EllLam19}), on the 135  instances of quadratically constrained quadratic programs from~\cite{BST09} called $unitbox$.  \\

\noindent \textit{Experimental environment}\\
 Our experiments were carried out on a server with $2$ CPU Intel Xeon each of them having $12$ cores and $2$ threads of $2.5$ GHz and $4*16$ GB of RAM using a Linux operating system. For all algorithms, we use the multi-threading version of \texttt{Cplex 12.7} with up to 48 threads. For methods \miqcr~and \miqcrt, we used the solver \csdp~(\cite{csdp}) together with the Conic Bundle library~(\cite{cb}) for solving semi-definite programs, as described in Section~\ref{sec:bundle}. We used the C interface of the solver \cplex(~\cite{cplex127}) for solving the quadratic convex relaxations at each node of the search tree. For computing feasible local solutions, we use the local solver \ipopt(~\cite{ipopt}). Parameter $p$ is set to $0.4\cdot\lvert \setMC \rvert$ for \miqcr, and to $0.04\cdot \lvert \setMC \cup \setGT \rvert$  for \miqcrt.\\

\noindent \textit{Results for the $unitbox$ instances}\\
Each instance  from~\cite{BST09} consists in minimizing a quadratic function of $n$ continuous variables in the interval $[0,1]$, subject to $m$ quadratic inequalities. For the considered instances, $n$ varies from 8 to 50, and $m$ from 8 to 100.  We set the time limit to 2 hours. In Figure~1, we present the performance profile of the CPU times for methods \miqcrt, \miqcr, and the solvers \glomiqo, and \baron~for the $unitbox$ instances. We observe that  \miqcrt~and \miqcr~outperform the compared solvers in terms of CPU time and number of instances solved. In fact, \texttt{BARON} solve 109 instances, \texttt{GloMIQO} solves 110 instances, \miqcr~solves 119 instances, and  \miqcrt~solves 128 instances out of 135 within the time limit. Several additional remarks are in order: the initial gap is smaller for \miqcrt~than for \miqcr, since we pass from $1.63\%$ to  $1.18\%$ on average on the 109 instances solved by both methods. Surprisingly, the CPU time for solving the semi-definite relaxation is divided by 2 on average for \miqcrt~in comparison to \miqcr, this is due to the sub-gradients considered in the course of the Conic Bundle algorithms that can be different for the two methods.  Another consequence of the use of the new inequalities within the \bb~process is the significant reduction of the number of nodes (by a factor $6.6$). Hence, we can also see a reduction of CPU time for this phase, where we pass from 390 seconds to 85 seconds on average. \\

\begin{figure}
\begin{center}
 \includegraphics[width=15cm]{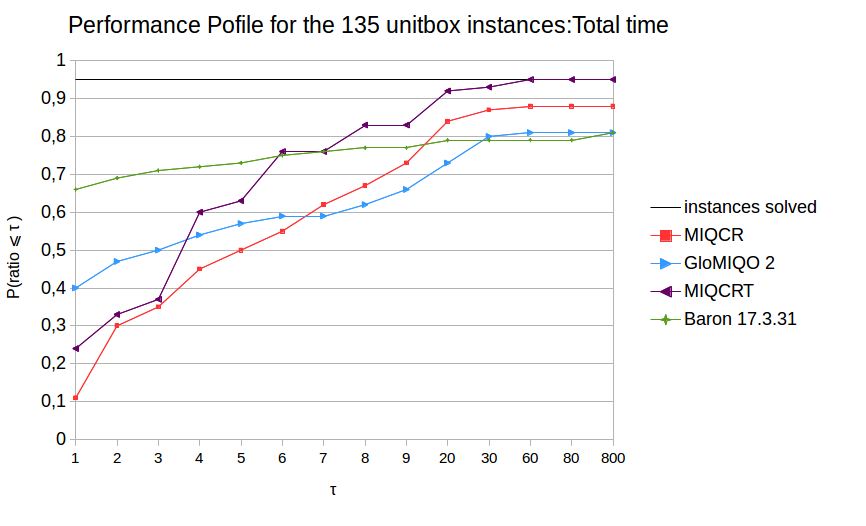}
\label{perfunitbox}
 \caption{Performance profile of the total time for the $unitbox$ instances with $n=8$ to $50$ with a time limit of 2 hours. }
\end{center}
\end{figure}

In Table~1, we present the detailed total CPU times for each method. Each line corresponds to one instance. We observe that \texttt{GloMIQO} or \texttt{BARON} are faster on most of the small and/or sparse instances, while \miqcr~and \miqcrt~are faster on large and dense instances.

\begin{landscape}
  
\begin{table}
\centering
\begin{scriptsize}
  \begin{tabular}{|l|c|c|c|c||l|c|c|c|c||l|c|c|c|c|} \hline
\textit{Name} &\miqcr& \miqcrt & \texttt{BARON}& \texttt{GloMIQO} &\textit{Name} &\miqcr& \miqcrt & \texttt{BARON}& \texttt{GloMIQO}&\textit{Name} &\miqcr& \miqcrt & \texttt{BARON}& \texttt{GloMIQO} \\ \hline
\textit{8\_12\_1\_25}&1&2&0&0&\textit{20\_40\_1\_25}&3&2&0&0&\textit{40\_60\_1\_25}&302&487&414&101\\ \hline
\textit{8\_12\_2\_25}&1&2&0&0&\textit{20\_40\_2\_25}&2&7&0&0&\textit{40\_60\_2\_25}&13&31&45&3\\ \hline
\textit{8\_12\_3\_25}&0&2&0&0&\textit{20\_40\_3\_25}&4&5&0&0&\textit{40\_60\_3\_25}&99&80&106&16\\ \hline
\textit{8\_16\_1\_25}&0&2&0&0&\textit{20\_40\_1\_50}&5&6&5&1&\textit{40\_60\_1\_50}&5516&190&282&1822\\ \hline
\textit{8\_16\_2\_25}&1&4&0&0&\textit{20\_40\_2\_50}&8&14&8&1&\textit{40\_60\_2\_50}&139&62&176&3616\\ \hline
\textit{8\_16\_3\_25}&0&2&0&0&\textit{20\_40\_3\_50}&1&3&0&0&\textit{40\_60\_3\_50}&3238&61&546&-\\ \hline
\textit{8\_8\_1\_25}&0&2&0&0&\textit{20\_40\_1\_100}&19&9&60&3&\textit{40\_60\_1\_100}&-&2913&-&-\\ \hline
\textit{8\_8\_2\_25}&1&3&0&0&\textit{20\_40\_2\_100}&18&12&50&2&\textit{40\_60\_2\_100}&3133&449&-&2182\\ \hline
\textit{8\_8\_3\_25}&0&2&0&0&\textit{20\_40\_3\_100}&3&4&1&0&\textit{40\_60\_3\_100}&7203&2663&-&-\\ \hline
\textit{10\_10\_1\_50}&0&2&0&0&\textit{28\_28\_1\_25}&2&3&0&0&\textit{40\_80\_1\_25}&445&208&1184&1558\\ \hline
\textit{10\_10\_2\_50}&0&2&0&0&\textit{28\_28\_2\_25}&2&7&0&0&\textit{40\_80\_2\_25}&326&292&338&107\\ \hline
\textit{10\_10\_3\_50}&1&3&0&0&\textit{28\_28\_3\_25}&7&24&2&1&\textit{40\_80\_3\_25}&1248&181&203&86\\ \hline
\textit{10\_10\_1\_100}&1&2&0&0&\textit{28\_42\_1\_25}&10&19&2&0&\textit{40\_80\_1\_50}&190&139&112&85\\ \hline
\textit{10\_10\_2\_100}&0&2&0&0&\textit{28\_42\_2\_25}&12&6&0&0&\textit{40\_80\_2\_50}&154&51&202&4785\\ \hline
\textit{10\_10\_3\_100}&1&3&0&0&\textit{28\_42\_3\_25}&1&4&0&0&\textit{40\_80\_3\_50}&4086&1014&-&-\\ \hline
\textit{10\_15\_1\_50}&1&2&0&0&\textit{28\_56\_1\_25}&4&8&1&1&\textit{40\_80\_1\_100}&-&250&6957&741\\ \hline
\textit{10\_15\_2\_50}&0&2&0&0&\textit{28\_56\_2\_25}&4&8&0&1&\textit{40\_80\_2\_100}&931&286&-&654\\ \hline
\textit{10\_15\_3\_50}&1&2&0&0&\textit{28\_56\_3\_25}&19&21&8&1&\textit{40\_80\_3\_100}&-&-&-&-\\ \hline
\textit{10\_15\_1\_100}&1&3&0&0&\textit{30\_30\_1\_50}&38&26&36&3&\textit{48\_48\_1\_25}&411&211&2140&1958\\ \hline
\textit{10\_15\_2\_100}&0&2&0&0&\textit{30\_30\_2\_50}&60&17&20&12&\textit{48\_48\_2\_25}&161&100&108&57\\ \hline
\textit{10\_15\_3\_100}&1&2&0&0&\textit{30\_30\_3\_50}&30&5&1&4&\textit{48\_48\_3\_25}&1400&671&393&452\\ \hline
\textit{10\_20\_1\_50}&0&3&0&0&\textit{30\_30\_1\_100}&813&269&-&202&\textit{48\_72\_1\_25}&64&88&85&32\\ \hline
\textit{10\_20\_2\_50}&3&2&0&0&\textit{30\_30\_2\_100}&38&13&9&3&\textit{48\_72\_2\_25}&322&126&206&37\\ \hline
\textit{10\_20\_3\_50}&0&3&0&0&\textit{30\_30\_3\_100}&49&8&175&51&\textit{48\_72\_3\_25}&546&151&112&180\\ \hline
\textit{10\_20\_1\_100}&2&2&0&0&\textit{30\_45\_1\_50}&106&25&24&7&\textit{48\_96\_1\_25}&268&186&381&33\\ \hline
\textit{10\_20\_2\_100}&2&2&0&0&\textit{30\_45\_2\_50}&49&18&11&3&\textit{48\_96\_2\_25}&4710&547&-&-\\ \hline
\textit{10\_20\_3\_100}&0&3&0&0&\textit{30\_45\_3\_50}&64&17&38&4&\textit{48\_96\_3\_25}&1060&575&1420&2879\\ \hline
\textit{20\_20\_1\_25}&1&4&0&0&\textit{30\_45\_1\_100}&184&42&210&26&\textit{50\_50\_1\_50}&-&315&-&-\\ \hline
\textit{20\_20\_2\_25}&0&3&0&0&\textit{30\_45\_2\_100}&375&85&360&61&\textit{50\_50\_2\_50}&-&2059&-&-\\ \hline
\textit{20\_20\_3\_25}&1&2&0&0&\textit{30\_45\_3\_100}&101&27&92&7&\textit{50\_50\_3\_50}&2593&862&-&-\\ \hline
\textit{20\_20\_1\_50}&1&3&0&0&\textit{30\_60\_1\_50}&833&103&188&323&\textit{50\_50\_1\_100}&2407&234&-&-\\ \hline
\textit{20\_20\_2\_50}&1&5&0&1&\textit{30\_60\_2\_50}&302&50&54&278&\textit{50\_50\_2\_100}&-&2680&-&-\\ \hline
\textit{20\_20\_3\_50}&1&2&0&0&\textit{30\_60\_3\_50}&57&15&9&7&\textit{50\_50\_3\_100}&-&3406&-&-\\ \hline
\textit{20\_20\_1\_100}&8&6&7&1&\textit{30\_60\_1\_100}&231&43&537&29&\textit{50\_75\_1\_50}&3255&1225&-&-\\ \hline
\textit{20\_20\_2\_100}&13&11&44&4&\textit{30\_60\_2\_100}&1602&230&5037&213&\textit{50\_75\_2\_50}&-&1233&-&-\\ \hline
\textit{20\_20\_3\_100}&3&3&0&0&\textit{30\_60\_3\_100}&272&34&894&26&\textit{50\_75\_3\_50}&5126&632&-&-\\ \hline
\textit{20\_30\_1\_25}&1&4&0&0&\textit{40\_40\_1\_25}&52&32&45&6&\textit{50\_75\_1\_100}&-&-&-&-\\ \hline
\textit{20\_30\_2\_25}&1&2&0&0&\textit{40\_40\_2\_25}&14&27&41&2&\textit{50\_75\_2\_100}&-&-&-&-\\ \hline
\textit{20\_30\_3\_25}&1&5&0&0&\textit{40\_40\_ 3\_25}&14&28&0&1&\textit{50\_75\_3\_100}&-&-&-&-\\ \hline
\textit{20\_30\_1\_50}&3&2&0&0&\textit{40\_40\_1\_50}&539&40&608&-&\textit{50\_100\_1\_50}&-&-&-&-\\ \hline
\textit{20\_30\_2\_50}&0&3&0&0&\textit{40\_40\_2\_50}&152&64&306&7034&\textit{50\_100\_2\_50}&823&324&-&-\\ \hline
\textit{20\_30\_3\_50}&2&3&0&0&\textit{40\_40\_3\_50}&574&94&2484&-&\textit{50\_100\_3\_50}&3337&1731&-&-\\ \hline
\textit{20\_30\_1\_100}&6&4&12&0&\textit{40\_40\_1\_100}&65&42&235&14&\textit{50\_100\_1\_100}&-&3808&-&-\\ \hline
\textit{20\_30\_2\_100}&10&7&25&2&\textit{40\_40\_2\_100}&297&72&2133&503&\textit{50\_100\_2\_100}&-&-&-&-\\ \hline
\textit{20\_30\_3\_100}&3&9&32&1&\textit{40\_40\_3\_100}&85&11&10&2&\textit{50\_100\_3\_100}&-&-&-&-\\ \hline
  \end{tabular}
  \end{scriptsize}
\label{tab:unit1}
\caption {Total CPU times of \miqcr, \miqcrt, \texttt{BARON~17.3.31} and \texttt{GloMIQO~2} for the 135 $unitbox$ instances. Time limit 2 hours.}

\end{table}
\end{landscape}

\section{Conclusion}
\label{sec:conc}
We consider  the general problem $\QP$ of minimizing a quadratic function  subject to quadratic constraints where the variables are continuous. In this paper, we introduce 12 General Triangle inequalities and we prove that they cut feasible solutions of the McCormick envelopes. In fact, these inequalities are an extension of the Triangle inequalities to the case of general lower and upper bounds on the variables. Then, we show how we can compute a quadratic convex relaxation which optimal value is equal to the "Shor's plus RLT plus Triangle" semi-definite relaxation $(SDP)$.
Since there is a huge number of these inequalities, we then focus on selecting the $p$ most violated ones. In particular, we separate them during the heuristic solution of $(SDP)$.  We report computational results on 135 instances. These results show that the method allows us to solve 128 instances out of 135 within a time limit of 2 hours.
From a general outlook, these new inequalities can be used in any \bb~process based on the relaxation of the constraints $Y=xx^T$. Indeed, since  the upper and lower bounds $\ell$ and $u$ are involved within the general Triangles,  the relaxation will be again tighten in the course of the algorithm.

\bigskip
\bibliographystyle{informs2014}
\bibliography{mybib}
     \end{document}